\documentclass[a4paper,11pt,twoside,fleqn]{article}
\usepackage{amsmath,amssymb,amsthm,latexsym}
\newtheorem{theorem}{Theorem}[section]
\newtheorem{proposition}[theorem]{Proposition}
\newtheorem{corollary}[theorem]{Corollary}
\newtheorem{lemma}[theorem]{Lemma}

\newcommand{\dist}{\text{\rm dist}}
\newcommand{\Id}{\text{\rm I}}

\begin{document}

\title{On the construction of solutions to the Yang-Mills equations in higher dimensions}
\author{Simon Brendle}
\date{August 13, 2003}

\maketitle

\section{Introduction}

Let $M$ be a Riemannian manifold of dimension $n \geq 4$. A a connection $A$ on a vector bundle 
over $M$ is a Yang-Mills connection if the curvature $F_A$ satisfies 
\[D_A^* F_A = 0.\] 
This is the Euler-Lagrange equation for the functional 
\[E[A] = \int_M |F_A|^2.\] 
Let $\{A_k\}$ be a sequence of Yang-Mills connections with uniformly bounded energy, i.e. 
\[\sup_k E[A_k] < \infty.\] 
Then the set 
\[S = \bigg \{ x \in M: \lim_{k \to \infty} r^{4-n} \int_{B_r(x)} |F_{A_k}|^2 \geq \varepsilon_0 
\quad \text{for all $r > 0$} \bigg \}\] 
is called the blow-up set of the sequence $\{A_k\}$. G. Tian \cite{Ti} proved that the blow-up 
set $S$ is closed and $H^{n-4}$-rectifiable. Moreover, the energy densities satisfy 
\[|F_{A_k}|^2 \, dvol \rightharpoonup |F_{A_\infty}|^2 \, dvol + 8\pi^2 \, \Theta \, dH^{n-4} |
_S\] 
as $k \to \infty$, where $A_\infty$ is the limiting connection defined on $M \setminus S$, 
$\Theta$ denotes the density function, and $H^{n-4}$ is the $(n-4)$-dimensional Hausdorff 
measure. Furthermore, if the limiting connection $A_\infty$ is admissible, then the generalized 
mean curvature of $S$ is equal to $0$ (see \cite{Ti,TT}). This result generalizes a theorem of 
K. Uhlenbeck in dimension $4$. \\

In this paper, we consider a smooth minimal submanifold $S$ of dimension $n - 4$. Our aim is to 
construct a sequence $\{A_k\}$ of smooth Yang-Mills connections whose blow-up set is equal to 
$S$. \\

In the first step, we construct a suitable family of approximate solutions. To this end, we 
assume that the normal bundle of $S$ can be endowed with a complex structure $J$ and a complex 
volume form $\omega$. Each approximate solution is described by a set $(v,\lambda,J,\omega)$, 
where $v$ is a section of the normal bundle of $S$, $\lambda$ is a positive function on $S$, and 
$(J,\omega)$ is a $SU(2)$-structure on $NS$. \\

For every point $x \in S$, the normal components of an approximate solution $A$ coincide with 
the basic instanton on the fibre $NS_x$, i.e. 
\begin{align*} 
A(e_1^\perp) &= \frac{-(y - \varepsilon v)_2 \, i - (y - \varepsilon v)_3 \, j - (y - 
\varepsilon v)_4 \, k}{\varepsilon^2 \lambda^2 + |y - \varepsilon v|^2} \\ 
A(e_2^\perp) &= \frac{(y - \varepsilon v)_1 \, i - (y - \varepsilon v)_4 \, j + (y - \varepsilon 
v)_3 \, k}{\varepsilon^2 \lambda^2 + |y - \varepsilon v|^2} \\ 
A(e_3^\perp) &= \frac{(y - \varepsilon v)_4 \, i + (y - \varepsilon v)_1 \, j - (y - \varepsilon 
v)_2 \, k}{\varepsilon^2 \lambda^2 + |y - \varepsilon v|^2} \\ 
A(e_4^\perp) &= \frac{-(y - \varepsilon v)_3 \, i + (y - \varepsilon v)_2 \, j + (y - 
\varepsilon v)_1 \, k}{\varepsilon^2 \lambda^2 + |y - \varepsilon v|^2}, 
\end{align*} 
where $\{e_\alpha^\perp: 1 \leq \alpha \leq 4\}$ denotes a $SU(2)$ basis for the fibre $NS_x$. 
\\

Our aim is to deform the connection $A$ to a nearby connection $\tilde{A} = A + a$ such that 
$\tilde{A}$ is a solution of the Yang-Mills equations. \\

In the following, we denote by $h$ the second fundamental form of the submanifold $S$, and by 
$R$ the Riemann curvature tensor of $M$. \\

\begin{theorem}
Suppose that $H^1(M) = 0$. Then, for each $\varepsilon > 0$, there exists a mapping $\Xi
_\varepsilon$ whch assigns to each set of glueing data $(v,\lambda,J,\omega) \in \mathcal{C}^{2,
\gamma}(S)$ a section of the vector bundle $NS \oplus \mathbb{R} \oplus \Lambda_+^2 NS$ of class 
$\mathcal{C}^\gamma(S)$ such that the following holds. \\

(i) If $(v,\lambda,J,\omega)$ is a set of glueing data such that 
\[\|v\|_{\mathcal{C}^{2,\gamma}(S)} \leq K,\] 
\[\|\lambda\|_{\mathcal{C}^{2,\gamma}(S)} \leq K, \qquad \inf \lambda \geq 1,\] 
\[\|(J,\omega)\|_{\mathcal{C}^{2,\gamma}(S)} \leq K,\] 
then we have the estimate 
\begin{align*} 
\bigg \| &\Xi_\varepsilon(v,\lambda,J,\omega) \\ 
&- \bigg ( \Delta v_\rho + \sum_{i,j=1}^{n-4} \sum_{\rho,\sigma=1}^4 h_{ij,\rho} \, h_{ij,
\sigma} \, v_\sigma + \sum_{i=1}^{n-4} \sum_{\rho,\sigma=1}^4 R_{i\rho\sigma i} \, v_\sigma, \\ 
&\hspace{7.5mm} \frac{1}{\lambda} \, \Delta \lambda + \frac{1}{4} \sum_{i,j=1}^{n-4} \sum_{\rho
=1}^4 h_{ij,\rho} \, h_{ij,\rho} + \frac{1}{4} \sum_{i=1}^{n-4} \sum_{\rho=1}^4 R_{i\rho\rho i} 
- \frac{1}{4} \, |\theta|^2, \\ 
&\hspace{7.5mm} \frac{1}{\lambda^2} \, \sum_{i=1}^{n-4} \nabla_i(\lambda^2 \, \theta_{i,\rho
\sigma}) \bigg ) \bigg \|_{\mathcal{C}^\gamma(S)} \leq C \, \varepsilon^{\frac{1}{64}}. 
\end{align*} 

(ii) If $\Xi_\varepsilon(v,\lambda,J,\omega) = 0$, then the approximate solution $A$ 
corresponding to $(v,\lambda,J,\omega)$ can be deformed to a nearby connection $\tilde{A}$ 
satisfying $D_{\tilde{A}}^* F_{\tilde{A}} = 0$.
\end{theorem}

\vspace{2mm}

In Section 2, we recall some results about the linearized operator on $\mathbb{R}^4$. In 
particular, the kernel of the linearized operator on $\mathbb{R}^4$ is isomorphic to $\mathbb{R}
^4 \oplus \mathbb{R}^4 \oplus \Lambda_+^2 \mathbb{R}^4$ (compare \cite{DK}). \\

In Section 3, we study the mapping properties of a model operator on the product manifold 
$\mathbb{R}^{n-4} \times \mathbb{R}^4$. \\

In Section 4, we construct a family of approximate solutions of the Yang-Mills equations. 
More precisely, given any set of glueing data $(v,\lambda,J,\omega)$ satisfying 
\[\|v\|_{\mathcal{C}^{2,\gamma}(S)} \leq K,\] 
\[\|\lambda\|_{\mathcal{C}^{2,\gamma}(S)} \leq K, \qquad \inf \lambda \geq 1,\] 
\[\|(J,\omega)\|_{\mathcal{C}^{2,\gamma}(S)} \leq K,\] 
we construct a connection $A$ such that 
\[\|D_A^* F_A\|_{\mathcal{C}_3^\gamma(M)} \leq C \, \varepsilon^2.\] 
Here, the weighted H\"older space $\mathcal{C}_\nu^\gamma(M)$ is defined as 
\begin{align*} 
\|u\|_{\mathcal{C}_\nu^\gamma(M)} 
&= \sup \, (\varepsilon + \dist(p,S))^\nu \, |u(p)| \\ 
&+ \sup_{\begin{smallmatrix} 4\dist(p_1,p_2) \leq \\ \varepsilon + \dist(p_1,S) + \dist(p_2,S) 
\end{smallmatrix}} \, (\varepsilon + \dist(p_1,S) + \dist(p_2,S))^{\nu+\gamma} \frac{|u(p_1) - 
u(p_2)|}{\dist(p_1,p_2)^\gamma}. 
\end{align*} 

In Section 5, we derive uniform estimates for the operator $\mathbb{L}_A = L_A + D_A D_A^*$. 
Here, $L_A$ is the linearization of the Yang-Mills equations at an approximate solution $A$. The 
additional term $D_A D_A^*$ must be included because $L_A$ is not an elliptic operator. \\

To derive uniform estimates independent of $\varepsilon$, we need to restrict the operator 
$\mathbb{L}_A$ to a subspace $\mathcal{E}_\nu^\gamma(M) \subset \mathcal{C}_\nu^\gamma(M)$. A 
$1$-form $a$ belongs to $\mathcal{E}_\nu^\gamma(M)$ if 
\[\int_{NS_x} \sum_{\alpha=1}^4 \langle a(e_\alpha^\perp),F_A(X,e_\alpha^\perp) \rangle = 0\] 
for all $x \in S$ and all vector fields of the form 
\[X = \varepsilon \, w_\rho \, e_\rho^\perp + \mu \, (y - \varepsilon v)_\rho \, e_\rho^\perp 
+ r_{\rho\sigma} \, (y - \varepsilon v)_\sigma \, e_\rho^\perp\] 
with $w \in NS_x$, $\mu \in \mathbb{R}$, and $r \in \Lambda_+^2 NS_x$. \\

In Section 6, we apply the contraction mapping principle to deform the approximate solution $A$ 
to a nearby connection $\tilde{A} = A + a$ such that 
\[(\Id - \mathbb{P}) (D_{\tilde{A}}^* F_{\tilde{A}} + D_{\tilde{A}} D_{\tilde{A}}^* a) = 0,\] 
where $(\Id - \mathbb{P})$ is the fibrewise projection from $\mathcal{C}_\nu^\gamma(M)$ to the 
subspace $\mathcal{E}_\nu^\gamma(M)$. In particular, if the balancing condition 
\[\mathbb{P} (D_{\tilde{A}}^* F_{\tilde{A}} + D_{\tilde{A}} D_{\tilde{A}}^* a) = 0\] 
is satisfied, then $\tilde{A}$ is a Yang-Mills connection. \\

In Section 7, we calculate the leading term in the asymptotic expansion of  
\[\mathbb{P} (D_{\tilde{A}}^* F_{\tilde{A}} + D_{\tilde{A}} D_{\tilde{A}}^* a) = 0.\] 
This concludes the proof of Theorem 1.1. \\

An example is discussed in Section 8. \\

A related balancing condition occurs in the work of R. Schoen and D. Pollack \cite{Po} 
on the constant scalar curvature equation in conformal geometry. In this way, an infinite 
dimensional problem is reduced to solving a finite dimensional balancing condition. The 
balancing condition ensures that the energy is stationary with respect to variations of the 
glueing data. For the constant scalar curvature equation, an asymptotic expansion for the energy 
of a multi-peak solution was calculated in A. Bahri and J. M. Coron \cite{BC}. Similar results 
have been proved for constant mean curvature hypersurfaces (see \cite{Ka1,Ka2,KMP}). \\

In the examples mentioned above, the blow-up set consists of isolated points and each 
approximate solution is characterized by a finite-dimensional set of glueing data. In our 
situation, the space of approximate solutions is infinite-dimensional. This leads to technical 
difficulties. A similar problem occurs in the work of F. Pacard and M. Ritor\'e \cite{PR1} on 
the gradient theory of phase transitions, and in our earlier work on the Ginzburg-Landau 
equations in higher dimensions. \\

The author is grateful to Professor Gang Tian for discussions, and to the referee for valuable 
comments. \\

\section{The kernel of the linearized operator on $\mathbb{R}^4$}

The basic instanton on $\mathbb{R}^4$ is given by the formula 
\begin{align*} 
B_1 &= \frac{-y_2 \, \mathfrak{i} - y_3 \, \mathfrak{j} - y_4 \, \mathfrak{k}}{\varepsilon^2 + 
|y|^2} \\ 
B_2 &= \frac{y_1 \, \mathfrak{i} - y_4 \, \mathfrak{j} + y_3 \, \mathfrak{k}}{\varepsilon^2 + 
|y|^2} \\ 
B_3 &= \frac{y_4 \, \mathfrak{i} + y_1 \, \mathfrak{j} - y_2 \, \mathfrak{k}}{\varepsilon^2 + 
|y|^2} \\ 
B_4 &= \frac{-y_3 \, \mathfrak{i} + y_2 \, \mathfrak{j} + y_1 \, \mathfrak{k}}{\varepsilon^2 + 
|y|^2}, 
\end{align*} 
where $\{\mathfrak{i},\mathfrak{j},\mathfrak{k}\}$ is the standard basis of $\mathfrak{su}(2)$, 
i.e. 
\begin{align*} 
&\mathfrak{i}(\frac{\partial}{\partial y_1}) = -\frac{\partial}{\partial y_2}, \quad 
\mathfrak{i}(\frac{\partial}{\partial y_2}) = \frac{\partial}{\partial y_1}, \quad 
\mathfrak{i}(\frac{\partial}{\partial y_3}) = \frac{\partial}{\partial y_4}, \quad 
\mathfrak{i}(\frac{\partial}{\partial y_4}) = -\frac{\partial}{\partial y_3}, \\ 
&\mathfrak{j}(\frac{\partial}{\partial y_1}) = -\frac{\partial}{\partial y_3}, \quad 
\mathfrak{j}(\frac{\partial}{\partial y_2}) = -\frac{\partial}{\partial y_4}, \quad 
\mathfrak{j}(\frac{\partial}{\partial y_3}) = \frac{\partial}{\partial y_1}, \quad 
\mathfrak{j}(\frac{\partial}{\partial y_4}) = \frac{\partial}{\partial y_2}, \\  
&\mathfrak{k}(\frac{\partial}{\partial y_1}) = -\frac{\partial}{\partial y_4}, \quad 
\mathfrak{k}(\frac{\partial}{\partial y_2}) = \frac{\partial}{\partial y_3}, \quad 
\mathfrak{k}(\frac{\partial}{\partial y_3}) = -\frac{\partial}{\partial y_2}, \quad   
\mathfrak{k}(\frac{\partial}{\partial y_4}) = \frac{\partial}{\partial y_1}.
\end{align*} 
Note that $[\mathfrak{i},\mathfrak{j}] = 2 \, \mathfrak{k}$, $[\mathfrak{j},\mathfrak{k}] = 2 \, 
\mathfrak{i}$, $[\mathfrak{k},\mathfrak{i}] = 2 \, \mathfrak{j}$. The curvature of $B$ satisfies 
\begin{align*} 
F_{B,12} = -F_{B,34} &= \frac{2\varepsilon^2 \, \mathfrak{i}}{(\varepsilon^2 + |y|^2)^2} \\ 
F_{B,13} = -F_{B,42} &= \frac{2\varepsilon^2 \, \mathfrak{j}}{(\varepsilon^2 + |y|^2)^2} \\ 
F_{B,14} = -F_{B,23} &= \frac{2\varepsilon^2 \, \mathfrak{k}}{(\varepsilon^2 + |y|^2)^2}. 
\end{align*}  
We first recall some well-known facts about this solution. In the first step, we construct a 
frame which is asymptotically parallel as $|y| \to \infty$. \\

\begin{lemma}
Let 
\[u = (\varepsilon^2 + |y|^2)^{-\frac{1}{2}} \, \big ( \mu \, y_\rho + r_{\rho\sigma} \, 
y_\sigma \big ) \, \frac{\partial}{\partial y_\rho}\] 
for some $\mu \in \mathbb{R}$ and $r \in \Lambda_+^2 \mathbb{R}^4$. Then 
\[D_{B,\alpha} u = \varepsilon^2 \, (\varepsilon^2 + |y|^2)^{-\frac{3}{2}} \, \big ( \mu \, 
\delta_{\rho\alpha} + r_{\rho\alpha} \big ) \, \frac{\partial}{\partial y_\rho}.\] 
\end{lemma}

\textbf{Proof.} 
We only consider the case $\mu = 1$, $r = 0$. By definition of $B$, we have 
\begin{align*} 
(\varepsilon^2 + |y|^2)^{\frac{3}{2}} \, D_{B,1} u 
&= (\varepsilon^2 + |y|^2) \, \frac{\partial}{\partial y_1} 
- y_1^2 \, \frac{\partial}{\partial y_1} - y_1 y_2 \, \frac{\partial}{\partial y_2} 
- y_1 y_3 \, \frac{\partial}{\partial y_3} - y_1 y_4 \, \frac{\partial}{\partial y_4} \\ 
&+ y_1 \Big ( y_2 \, \frac{\partial}{\partial y_2} 
+ y_3 \, \frac{\partial}{\partial y_3} 
+ y_4 \, \frac{\partial}{\partial y_4} \Big ) \\ 
&+ y_2 \Big ( -y_2 \, \frac{\partial}{\partial y_1} 
+ y_3 \, \frac{\partial}{\partial y_4} 
- y_4 \, \frac{\partial}{\partial y_3} \Big ) \\ 
&+ y_3 \Big ( -y_2 \, \frac{\partial}{\partial y_4} 
- y_3 \, \frac{\partial}{\partial y_1} 
+ y_4 \, \frac{\partial}{\partial y_2} \Big ) \\ 
&+ y_4 \Big ( y_2 \, \frac{\partial}{\partial y_3} 
- y_3 \, \frac{\partial}{\partial y_2} 
- y_4 \, \frac{\partial}{\partial y_1} \Big ) \\ 
&= \varepsilon^2 \, \frac{\partial}{\partial y_1}, 
\end{align*} 
\begin{align*} 
(\varepsilon^2 + |y|^2)^{\frac{3}{2}} \, D_{B,2} u 
&= (\varepsilon^2 + |y|^2) \, \frac{\partial}{\partial y_2} 
- y_2 y_1 \, \frac{\partial}{\partial y_1} - y_2^2 \, \frac{\partial}{\partial y_2} 
- y_2 y_3 \, \frac{\partial}{\partial y_3} - y_2 y_4 \, \frac{\partial}{\partial y_4} \\ 
&+ y_1 \Big ( -y_1 \, \frac{\partial}{\partial y_2} 
+ y_4 \, \frac{\partial}{\partial y_3} 
- y_3 \, \frac{\partial}{\partial y_4} \Big ) \\ 
&+ y_2 \Big ( y_1 \, \frac{\partial}{\partial y_1} 
+ y_4 \, \frac{\partial}{\partial y_4} 
+ y_3 \, \frac{\partial}{\partial y_3} \Big ) \\ 
&+ y_3 \Big ( y_1 \, \frac{\partial}{\partial y_4} 
- y_4 \, \frac{\partial}{\partial y_1} 
- y_3 \, \frac{\partial}{\partial y_2} \Big ) \\ 
&+ y_4 \Big ( -y_1 \, \frac{\partial}{\partial y_3} 
- y_4 \, \frac{\partial}{\partial y_2} 
+ y_3 \, \frac{\partial}{\partial y_1} \Big ) \\ 
&= \varepsilon^2 \, \frac{\partial}{\partial y_2}, 
\end{align*} 
\begin{align*} 
(\varepsilon^2 + |y|^2)^{\frac{3}{2}} \, D_{B,3} u 
&= (\varepsilon^2 + |y|^2) \, \frac{\partial}{\partial y_3} 
- y_3 y_1 \, \frac{\partial}{\partial y_1} - y_3 y_2 \, \frac{\partial}{\partial y_2} 
- y_3^2 \, \frac{\partial}{\partial y_3} - y_3 y_4 \, \frac{\partial}{\partial y_4} \\ 
&+ y_1 \Big ( -y_4 \, \frac{\partial}{\partial y_2} 
- y_1 \, \frac{\partial}{\partial y_3} 
+ y_2 \, \frac{\partial}{\partial y_4} \Big ) \\ 
&+ y_2 \Big ( y_4 \, \frac{\partial}{\partial y_1} 
- y_1 \, \frac{\partial}{\partial y_4} 
- y_2 \, \frac{\partial}{\partial y_3} \Big ) \\ 
&+ y_3 \Big ( y_4 \, \frac{\partial}{\partial y_4} 
+ y_1 \, \frac{\partial}{\partial y_1} 
+ y_2 \, \frac{\partial}{\partial y_2} \Big ) \\ 
&+ y_4 \Big ( -y_4 \, \frac{\partial}{\partial y_3} 
+ y_1 \, \frac{\partial}{\partial y_2} 
- y_2 \, \frac{\partial}{\partial y_1} \Big ) \\ 
&= \varepsilon^2 \, \frac{\partial}{\partial y_3}, 
\end{align*} 
\begin{align*} 
(\varepsilon^2 + |y|^2)^{\frac{3}{2}} \, D_{B,4} u 
&= (\varepsilon^2 + |y|^2) \, \frac{\partial}{\partial y_4} 
- y_4 y_1 \, \frac{\partial}{\partial y_1} - y_4 y_2 \, \frac{\partial}{\partial y_2} 
- y_4 y_3 \, \frac{\partial}{\partial y_3} - y_4^2 \, \frac{\partial}{\partial y_4} \\ 
&+ y_1 \Big ( y_3 \, \frac{\partial}{\partial y_2} 
- y_2 \, \frac{\partial}{\partial y_3} 
- y_1 \, \frac{\partial}{\partial y_4} \Big ) \\ 
&+ y_2 \Big ( -y_3 \, \frac{\partial}{\partial y_1} 
- y_2 \, \frac{\partial}{\partial y_4} 
+ y_1 \, \frac{\partial}{\partial y_3} \Big ) \\ 
&+ y_3 \Big ( -y_3 \, \frac{\partial}{\partial y_4} 
+ y_2 \, \frac{\partial}{\partial y_1} 
- y_1 \, \frac{\partial}{\partial y_2} \Big ) \\ 
&+ y_4 \Big ( y_3 \, \frac{\partial}{\partial y_3} 
+ y_2 \, \frac{\partial}{\partial y_2} 
+ y_1 \, \frac{\partial}{\partial y_1} \Big ) \\ 
&= \varepsilon^2 \, \frac{\partial}{\partial y_4}. 
\end{align*} 
The remaining cases are left to the reader. \\

\begin{lemma}
For every $r \in \Lambda_+^2 \mathbb{R}^4$, we have the identity 
\[D_B(r_{\rho\sigma} \, y_\sigma \, B_\rho) = -F_B \Big ( r_{\rho \sigma} \, y_\sigma \, \frac
{\partial}{\partial y_\rho},\cdot \Big ).\] 
\end{lemma}

\textbf{Proof.}
By direct calculation, one can see that 
\[r_{\rho\sigma} \, y_\sigma \, \partial_\rho B_\alpha + r_{\rho\alpha} \, B_\rho = 0\] 
for all $r \in \Lambda_+^2 \mathbb{R}^4$. This implies 
\begin{align*} 
r_{\rho\sigma} \, y_\sigma \, F_{B,\rho\alpha} 
&= r_{\rho\sigma} \, y_\sigma \, (\partial_\rho B_\alpha - \partial_\alpha B_\rho + [B_\rho,B
_\alpha]) \\ 
&= r_{\rho\sigma} \, y_\sigma \, (\partial_\rho B_\alpha - D_{B,\alpha} B_\rho) \\ 
&= r_{\rho\sigma} \, y_\sigma \, \partial_\rho B_\alpha + r_{\rho\alpha} B_\rho - D_{B,\alpha} 
(r_{\rho\sigma} \, y_\sigma \, B_\rho) \\ 
&= -D_{B,\alpha} (r_{\rho\sigma} \, y_\sigma \, B_\rho). 
\end{align*} 
This proves the assertion. \\

In the second step, we study the linearized operator (cf. \cite{Ta1}), which we denote by $L_B$. 
The operator $L_B$ is given by the formula 
\[L_B a = D_B^* D_B a - *[F_B,a].\] 
Using the Weitzenb\"ock formula 
\[D_B^* D_B a + D_B D_B^* a = \nabla_B^* \nabla_B a - *[F_B,a],\] 
we obtain 
\[L_B a + D_B D_B^* a = \nabla_B^* \nabla_B a - 2 \, *[F_B,a].\] 
In particular, $L_B + D_B D_B^*$ is an elliptic operator. \\

\begin{lemma}
For every $a \in \Omega^1(\mathbb{R}^4)$ we have 
\[L_B a + D_B D_B^* a = 2 \, D_B^* P_+ D_B a + D_B D_B^* a.\] 
Hence, the operator $L_B + D_B D_B^*$ is positive semidefinite.
\end{lemma}

\textbf{Proof.} 
By definition of $L_B$, we have 
\begin{align*} 
L_B a 
&= D_B^* D_B a - *[F_B,a] \\ 
&= D_B^* D_B a - * D_B D_B a \\ 
&= D_B^* D_B a + D_B^* * D_B a \\ 
&= 2 \, D_B^* P_+ D_B a. 
\end{align*} 
From this the assertion follows. \\

Our aim is to describe the kernel of the operator $L_B + D_B D_B^*: \mathcal{C}_{1+\nu}^{2,
\gamma}(\mathbb{R}^4) \to \mathcal{C}_{3+\nu}^\gamma(\mathbb{R}^4)$. \\

\begin{proposition}
Let 
\[a = F_B \Big ( \varepsilon \, w_\rho \, \frac{\partial}{\partial y_\rho} + \mu \, y_\rho \, 
\frac{\partial}{\partial y_\rho} + r_{\rho\sigma} \, y_\sigma \, \frac{\partial}{\partial 
y_\rho},\cdot \Big )\] 
for some $w \in \mathbb{R}^4$, $\mu \in \mathbb{R}$, and $r \in \Lambda_+^2 \mathbb{R}^4$. Then 
$a$ satisfies $P_+ D_B a = 0$ and $D_B^* a = 0$.
\end{proposition}

\textbf{Proof.}
Using the Bianchi identity, we obtain 
\[D_B a = \varepsilon \, w_\rho \, D_{B,\rho} F_B + \mu \, y_\rho \, D_{B,\rho} F_B + r_{\rho
\sigma} \, y_\sigma \, D_{B,\rho} F_B + 2\mu \, F_B.\] 
Since $F_B \in \Lambda_-^2 \mathbb{R}^4$, this implies $P_+ D_B a = 0$. Similarly, we obtain 
\[D_B^* a = (D_B^* F_B) \Big ( \varepsilon \, w_\rho \, \frac{\partial}{\partial y_\rho} + \mu 
\, y_\rho \, \frac{\partial}{\partial y_\rho} + r_{\rho\sigma} \, y_\sigma \, \frac
{\partial}{\partial y_\rho} \Big ) = 0.\] 
This proves the assertion. \\

\begin{corollary}
Let \[a = F_B \Big ( \varepsilon \, w_\rho \, \frac{\partial}{\partial y_\rho} + \mu \, y_\rho 
\, \frac{\partial}{\partial y_\rho} + r_{\rho\sigma} \, y_\sigma \, \frac{\partial}{\partial 
y_\rho},\cdot \Big )\] 
for some $w \in \mathbb{R}^4$, $\mu \in \mathbb{R}$, and $r \in \Lambda_+^2 \mathbb{R}^4$. Then 
$a$ satisfies $L_B a + D_B D_B^* a = 0$.
\end{corollary}

\vspace{2mm}

We now prove the converse statement. Due to the conformal invariance of the Yang-Mills equations 
in dimension $4$, we may lift the problem on $S^4$. The kernel of the linearized operator on 
$S^4$ is described in the following result which is well-known. \\

\begin{lemma}
Let $a$ be a $1$-form on $S^4$ such that $L_{B,g_{S^4}} a = 0$. Then there exists some $w \in 
\mathbb{R}^4$, $\mu \in \mathbb{R}$, and an infinitesimal gauge transformation $u$ such that 
\[a = F_B \Big ( \varepsilon \, w_\rho \, \frac{\partial}{\partial y_\rho} + \mu \, y_\rho \, 
\frac{\partial}{\partial y_\rho},\cdot \Big ) + D_B u.\]
\end{lemma}

\vspace{2mm}

\begin{proposition}
Let $0 < \nu < 1$. Assume that $a \in \mathcal{C}_{1+\nu}^{2,\gamma}(\mathbb{R}^4)$ satisfies 
\[L_B a + D_B D_B^* a = 0.\] 
Then $a$ is of the form 
\[a = F_B \Big ( \varepsilon \, w_\rho \, \frac{\partial}{\partial y_\rho} + \mu \, y_\rho \, 
\frac{\partial}{\partial y_\rho} + r_{\rho\sigma} \, y_\sigma \, \frac{\partial}{\partial 
y_\rho},\cdot \Big )\] 
for some $w \in \mathbb{R}^4$, $\mu \in \mathbb{R}$, and $r \in \Lambda_+^2 \mathbb{R}^4$.
\end{proposition}

\textbf{Proof.} 
Using integration by parts, we deduce that $L_B a = 0$ and $D_B^* a = 0$. We now lift the 
problem to $S^4$. The round metric on $S^4$ is given by 
\[g_{S^4} = \frac{4\varepsilon^2}{(\varepsilon^2 + |y|^2)^2} \, g_{\mathbb{R}^4}.\] 
Using the estimate 
\[|a(y)|_{\mathbb{R}^4} = O(|y|^{-1-\nu}),\] 
we obtain 
\[|a(y)|_{S^4} = \frac{\varepsilon^2 + |y|^2}{2\varepsilon} \, |a(y)|_{\mathbb{R}^4} = O(|y|^{1-
\nu}).\] 
The conformal invariance of the Yang-Mills equations implies that $L_{B,g_{S^4}} a = 0$. From 
this it follows that 
\[a = F_B \Big ( \varepsilon \, w_\rho \, \frac{\partial}{\partial y_\rho} + \mu \, y_\rho \, 
\frac{\partial}{\partial y_\rho},\cdot \Big ) + D_B u,\] 
where $w \in \mathbb{R}^4$, $\mu \in \mathbb{R}$, and $u$ denotes an infinitesimal gauge 
transformation. The function $u$ can be written as 
\[u = -r_{\rho\sigma} \, y_\sigma \, B_\rho + u_0,\] 
where $r \in \Lambda_+^2 \mathbb{R}^4$ and $u_0 = O(|y|^{-\nu})$. Therefore, we obtain 
\[a = F_B \Big ( \varepsilon \, w_\rho \, \frac{\partial}{\partial y_\rho} + \mu \, y_\rho \, 
\frac{\partial}{\partial y_\rho} + r_{\rho\sigma} \, y_\sigma \, \frac{\partial}{\partial 
y_\rho},\cdot \Big ) + D_B u_0.\] 
Since $D_B^* a = 0$, it follows that $D_B^* D_B u_0 = 0$. Since $u_0 \in \mathcal{C}_\nu^{3,
\gamma}(\mathbb{R}^4)$, we conclude that $u_0 = 0$. Therefore, we obtain 
\[a = F_B \Big ( \varepsilon \, w_\rho \, \frac{\partial}{\partial y_\rho} + \mu \, y_\rho \, 
\frac{\partial}{\partial y_\rho} + r_{\rho\sigma} \, y_\sigma \, \frac{\partial}{\partial 
y_\rho},\cdot \Big ).\] 
This proves the assertion. \\

\section{The model problem on $\mathbb{R}^{n-4} \times \mathbb{R}^4$}

Let $B$ be a connection on $\mathbb{R}^{n-4} \times \mathbb{R}^4$ which is invariant under 
translations along the $\mathbb{R}^{n-4}$ factor and agrees with the basic instanton along the 
$\mathbb{R}^4$ factor. This implies  
\begin{align*} 
B(e_1^\perp) &= \frac{-y_2 \, \mathfrak{i} - y_3 \, \mathfrak{j} - y_4 \, \mathfrak
{k}}{\varepsilon^2 + |y|^2} \\ 
B(e_2^\perp) &= \frac{y_1 \, \mathfrak{i} - y_4 \, \mathfrak{j} + y_3 \, \mathfrak
{k}}{\varepsilon^2 + |y|^2} \\ 
B(e_3^\perp) &= \frac{y_4 \, \mathfrak{i} + y_1 \, \mathfrak{j} - y_2 \, \mathfrak
{k}}{\varepsilon^2 + |y|^2} \\ 
B(e_4^\perp) &= \frac{-y_3 \, \mathfrak{i} + y_2 \, \mathfrak{j} + y_1 \, \mathfrak
{k}}{\varepsilon^2 + |y|^2}, 
\end{align*} 
where 
\begin{align*} 
&\mathfrak{i}(e_1^\perp) = -e_2^\perp, \quad \mathfrak{i}(e_2^\perp) = e_1^\perp, \quad 
\mathfrak{i}(e_3^\perp) = e_4^\perp, \quad \mathfrak{i}(e_4^\perp) = -e_3^\perp, \\ 
&\mathfrak{j}(e_1^\perp) = -e_3^\perp, \quad \mathfrak{j}(e_2^\perp) = -e_4^\perp, \quad 
\mathfrak{j}(e_3^\perp) = e_1^\perp, \quad \mathfrak{j}(e_4^\perp) = e_2^\perp, \\ 
&\mathfrak{k}(e_1^\perp) = -e_4^\perp, \quad \mathfrak{k}(e_2^\perp) = e_3^\perp, \quad 
\mathfrak{k}(e_3^\perp) = -e_2^\perp, \quad \mathfrak{k}(e_4^\perp) = e_1^\perp. 
\end{align*} 
Furthermore, $B(e_i) = 0$ for $1 \leq i \leq 4$. \\

The linearized operator satisfies 
\[L_B a = D_B^* D_B a + (-1)^n \, *[*F_B,a]\] 
Using the Weitzenb\"ock formula 
\[D_B^* D_B a + D_B D_B^* a = \nabla_B^* \nabla_B a + (-1)^n \, *[*F_B,a],\] 
we obtain 
\[L_B a + D_B D_B^* a = \nabla_B^* \nabla_B a + (-1)^n \, 2 \, *[*F_B,a],\] 
For abbreviation, let $\mathbb{L}_B = L_B + D_B D_B^*$. \\

We define the weighted H\"older space $\mathcal{C}_\nu^\gamma(\mathbb{R}^{n-4} \times \mathbb{R}
^4)$ by 
\begin{align*} 
\|u\|_{\mathcal{C}_\nu^\gamma(\mathbb{R}^{n-4} \times \mathbb{R}^4)} 
&= \sup \, (\varepsilon + |y|)^\nu \, |u(x,y)| \\ 
&+ \sup_{\begin{smallmatrix} 4 (|x_1 - x_2| + |y_1 - y_2|) \leq \\ \varepsilon + |y_1| + |y_2| 
\end{smallmatrix}} \, (\varepsilon + |y_1| + |y_2|)^{\nu+\gamma} \frac{|u(x_1,y_1) - u(x_2,y
_2)|}{(|x_1 - x_2| + |y_1 - y_2|)^\gamma}. 
\end{align*} 
More generally, we define 
\[\|u\|_{\mathcal{C}_\nu^{k,\gamma}(\mathbb{R}^{n-4} \times \mathbb{R}^4)} = \sum_{l=0}^k \|
\nabla^l u\|_{\mathcal{C}_{\nu+l}^\gamma(\mathbb{R}^{n-4} \times \mathbb{R}^4)}.\] 

Let $\mathcal{E}_\nu^{k,\gamma}(\mathbb{R}^{n-4} \times \mathbb{R}^4)$ be the set of all $a \in 
\Omega^1(\mathbb{R}^{n-4} \times \mathbb{R}^4)$ such that $a \in \mathcal{C}_\nu^{k,\gamma}
(\mathbb{R}^{n-4} \times \mathbb{R}^4)$ and 
\[\int_{\{x\} \times \mathbb{R}^4} \sum_{\alpha=1}^4 \langle a(e_\alpha^\perp),F_B(X,e_\alpha
^\perp) \rangle = 0\] 
for all $x \in \mathbb{R}^{n-4}$ and all vector fields of the form 
\[X = \varepsilon \, w_\rho \, e_\rho^\perp + \mu \, y_\rho \, e_\rho^\perp + r_{\rho\sigma} \, 
y_\sigma \, e_\rho^\perp\] 
with $w \in \mathbb{R}^4$, $\mu \in \mathbb{R}$, and $r \in \Lambda_+^2 \mathbb{R}^4$. \\

\begin{proposition}
The operator $\mathbb{L}_B$ maps $\mathcal{E}_{1+\nu}^{2,\gamma}(\mathbb{R}^{n-4} \times \mathbb
{R}^4)$ into $\mathcal{E}_{3+\nu}^\gamma(\mathbb{R}^{n-4} \times \mathbb{R}^4)$.
\end{proposition}

\textbf{Proof.} 
It is obvious from the definition that $\mathbb{L}_B$ maps $\mathcal{C}_{1+\nu}^{2,\gamma}
(\mathbb{R}^{n-4} \times \mathbb{R}^4)$ into $\mathcal{C}_{3+\nu}^\gamma(\mathbb{R}^{n-4} \times 
\mathbb{R}^4)$. We now assume that $a \in \mathcal{C}_{1+\nu}^{2,\gamma}(\mathbb{R}^{n-4} \times 
\mathbb{R}^4)$ satisfies 
\[\int_{\{x\} \times \mathbb{R}^4} \sum_{\alpha=1}^4 \langle a(e_\alpha^\perp),F_B(X,e_\alpha
^\perp) \rangle = 0\] 
for all $x \in \mathbb{R}^{n-4}$ and all vector fields of the form 
\[X = \varepsilon \, w_\rho \, e_\rho^\perp + \mu \, y_\rho \, e_\rho^\perp + r_{\rho\sigma} \, 
y_\sigma \, e_\rho^\perp\] 
with $w \in \mathbb{R}^4$, $\mu \in \mathbb{R}$, and $r \in \Lambda_+^2 \mathbb{R}^4$. Taking 
derivatives in horizontal direction, we obtain 
\[\int_{\{x\} \times \mathbb{R}^4} \sum_{\alpha=1}^4 \sum_{j=1}^{n-4} \langle \partial_j 
\partial_j a(e_\alpha^\perp),F_B(X,e_\alpha^\perp) \rangle = 0.\] 
Furthermore, integration by parts gives 
\begin{align*} 
&\int_{\{x\} \times \mathbb{R}^4} \sum_{\alpha=1}^4 \sum_{\beta=1}^4 \langle D_{B,e_\beta^\perp} 
D_{B,e_\beta^\perp} a(e_\alpha^\perp),F_B(X,e_\alpha^\perp) \rangle \\ 
&+ 2 \int_{\{x\} \times \mathbb{R}^4} \sum_{\alpha=1}^4 \sum_{\beta=1}^4 \langle [F_B(e_\alpha
^\perp,e_\beta^\perp),a(e_\beta^\perp)],F_B(X,e_\alpha^\perp) \rangle \\ 
&= \int_{\{x\} \times \mathbb{R}^4} \sum_{\alpha=1}^4 \sum_{\beta=1}^4 \langle a(e_\alpha
^\perp),D_{B,e_\beta^\perp} D_{B,e_\beta^\perp} F_B(X,e_\alpha^\perp) \rangle \\ 
&+ 2 \int_{\{x\} \times \mathbb{R}^4} \sum_{\alpha=1}^4 \sum_{\beta=1}^4 \langle a(e_\alpha
^\perp),[F_B(e_\alpha^\perp,e_\beta^\perp),F_B(X,e_\alpha^\perp)] \rangle \\ 
&= 0. 
\end{align*} 
Thus, we conclude that 
\[\int_{\{x\} \times \mathbb{R}^4} \sum_{\alpha=1}^4 \langle (\mathbb{L}_B a)(e_\alpha^\perp),F
_B(X,e_\alpha^\perp) \rangle = 0\] 
for all $x \in \mathbb{R}^{n-4}$ and all vector fields of the form 
\[X = \varepsilon \, w_\rho \, e_\rho^\perp + \mu \, y_\rho \, e_\rho^\perp + r_{\rho\sigma} \, 
y_\sigma \, e_\rho^\perp\] 
with $w \in \mathbb{R}^4$, $\mu \in \mathbb{R}$, and $r \in \Lambda_+^2 \mathbb{R}^4$. \\

\begin{proposition}
Let $0 < \nu < 1$, $b \in \mathcal{C}_{3+\nu}^\gamma(\mathbb{R}^4)$, and $\eta \in \mathcal{S}
(\mathbb{R}^{n-4})$. Moreover, assume that the Fourier transform of $\eta$ satisfies $\hat{\eta}
(\xi) = 0$ for $|\xi| \leq \delta$ for some $\delta > 0$. Then there exists a $1$-form $a \in 
\mathcal{C}_{1+\nu}^{2,\gamma}(\mathbb{R}^{n-4} \times \mathbb{R}^4)$ such that 
\[\mathbb{L}_B a = \eta(x) \, b(y).\]
\end{proposition}

\textbf{Proof.} 
We perform a Fourier transformation in the $\mathbb{R}^{n-4}$ variables. Let 
\[\eta(x) = \int_{\mathbb{R}^{n-4}} e^{ix\xi} \, \hat{\eta}(\xi) \, d\xi.\] 
For every $\xi \in \mathbb{R}^{n-4}$, there exists a $1$-form $\hat{a}(\xi,\cdot) \in \mathcal
{C}_{1+\nu}^{2,\gamma}(\mathbb{R}^4)$ such that 
\[\sum_{\beta=1}^4 D_{B,\beta} D_{B,\beta} \hat{a}_\alpha(\xi,y) + 2 \sum_{\beta=1}^4 [F_{B,
\alpha\beta},\hat{a}_\beta(\xi,y)] - |\xi|^2 \, \hat{a}_\alpha(\xi,y) = -b_\alpha(y)\] 
for $1 \leq \alpha \leq 4$ and 
\[\sum_{\beta=1}^4 D_{B,\beta} D_{B,\beta} \hat{a}_i(\xi,y) - |\xi|^2 \, \hat{a}_i(\xi,y) = -b_i
(y)\] 
for $1 \leq i \leq n - 4$. We now define a $1$-form $a$ by 
\[a_\alpha(x,y) = \int_{\mathbb{R}^{n-4}} e^{ix\xi} \, \hat{\eta}(\xi) \, \hat{a}_\alpha(\xi,y) 
\, d\xi\] 
for $1 \leq \alpha \leq 4$ and 
\[a_i(x,y) = \int_{\mathbb{R}^{n-4}} e^{ix\xi} \, \hat{\eta}(\xi) \, \hat{a}_i(\xi,y) \, d\xi\] 
for $1 \leq i \leq n-4$. Then the $1$-form $a$ satisfies 
\[\sum_{j=1}^{n-4} \partial_j \partial_j a_\alpha + \sum_{\beta=1}^4 D_{B,\beta} D_{B,\beta} a
_\alpha + 2 \sum_{\beta=1}^4 [F_{B,\alpha\beta},a_\beta] = -\eta(x) \, b_\alpha(y)\] 
for $1 \leq \alpha \leq 4$ and 
\[\sum_{j=1}^{n-4} \partial_j \partial_j a_i + \sum_{\beta=1}^4 D_{B,\beta} D_{B,\beta} a_i = 
-\eta(x) \, b_i(y)\] 
for $1 \leq i \leq n - 4$. From this we deduce that $\mathbb{L}_B a = \eta(x) \, b(y)$. \\

\begin{proposition}
Let $0 < \nu < 1$, and suppose that $a \in \mathcal{E}_{1+\nu}^{2,\gamma}(\mathbb{R}^{n-4} 
\times \mathbb{R}^4)$ satisfies $\mathbb{L}_B a = 0$. Then $a = 0$.
\end{proposition}

\textbf{Proof.} 
Let $b \in \mathcal{C}_{3+\nu}^\gamma(\mathbb{R}^4)$ and $\zeta \in \mathcal{S}(\mathbb{R}^{n-
4})$ be given. We define a function $\eta \in \mathcal{S}(\mathbb{R}^{n-4})$ by $\eta(x) = \zeta
(x + x_0) - \zeta(x)$. Then the Fourier transform of $\eta$ satisfies $\hat{\eta}(0) = 0$. We 
approximate $\eta$ by functions $\eta_\delta$ such that 
\[\hat{\eta}_\delta(\xi) = 0\] 
for $|\xi| \leq \delta$ and 
\[\hat{\eta}_\delta(\xi) = \hat{\eta}(\xi)\] 
for $|\xi| \geq 2\delta$. Using the condition $\hat{\eta}(0) = 0$, we obtain 
\[\|D^{n-4}(\hat{\eta} - \hat{\eta}_\delta)\|_{L^{\frac{p}{p-1}}(\mathbb{R}^{n-4})} \leq C \, 
\delta^{1-\frac{n-4}{p}},\] 
hence 
\[\big \| (1 + |x|)^{n-4} \, (\eta - \eta_\delta) \big \|_{L^p(\mathbb{R}^{n-4})} \leq C \, 
\delta^{1-\frac{n-4}{p}}\] 
for all $p \geq 2$. From this it follows that 
\begin{align*} 
\|\eta - \eta_\delta\|_{L^1(\mathbb{R}^{n-4})} 
&\leq \big \| (1 + |x|)^{-(n-4)} \big \|_{L^{\frac{p}{p-1}}(\mathbb{R}^{n-4})} \, \big \| (1 + 
|x|)^{n-4} \, (\eta - \eta_\delta) \big \|_{L^p(\mathbb{R}^{n-4})} \\ 
&\leq C \, \delta^{1-\frac{n-4}{p}} 
\end{align*} 
for all $p \geq 2$. This implies 
\[\|\eta - \eta_\delta\|_{L^1(\mathbb{R}^{n-4})} \to 0\] 
as $\delta \to 0$. \\

For each $\delta > 0$, the $1$-form $\eta_\delta(x) \, b(y)$ lies in the image of $\mathbb{L}
_B$. Since $a$ belongs to the kernel of $\mathbb{L}_B$, we obtain 
\[\int_{\mathbb{R}^{n-4} \times \mathbb{R}^4} \langle a(x,y),\eta_\delta(x) \, b(y) \rangle = 0.
\] 
Letting $\delta \to 0$, we obtain 
\[\int_{\mathbb{R}^{n-4} \times \mathbb{R}^4} \langle a(x,y),\eta(x) \, b(y) \rangle = 0,\] 
hence 
\[\int_{\mathbb{R}^{n-4} \times \mathbb{R}^4} \langle a(x,y),\zeta(x) \, b(y) \rangle = \int
_{\mathbb{R}^{n-4} \times \mathbb{R}^4} \langle a(x-x_0,y),\zeta(x) \, b(y) \rangle.\] 
Since $b$ and $\zeta$ are arbitrary, we conclude that $a(x,y) = a(x-x_0,y)$. Therefore, $a(x,y)$ 
is constant in $x$. Using Proposition 2.7, we obtain 
\[a = F_B(X,\cdot),\] 
where $X$ is a vector field of the form 
\[X = \varepsilon \, w_\rho \, e_\rho^\perp + \mu \, y_\rho \, e_\rho^\perp + r_{\rho\sigma} \, 
y_\sigma \, e_\rho^\perp = 0\] 
for suitable $w \in \mathbb{R}^4$, $\mu \in \mathbb{R}$, and $r \in \Lambda_+^2 \mathbb{R}^4$. 
This proves the assertion. \\

\begin{proposition}
Let $0 < \nu < 1$. Then we have the estimate 
\[\|a\|_{\mathcal{C}_{1+\nu}^{2,\gamma}(\mathbb{R}^{n-4} \times \mathbb{R}^4)} \leq C \, 
\|\mathbb{L}_B a\|_{\mathcal{C}_{3+\nu}^\gamma(\mathbb{R}^{n-4} \times \mathbb{R}^4)}\] 
for all $a \in \mathcal{E}_{1+\nu}^{2,\gamma}(\mathbb{R}^{n-4} \times \mathbb{R}^4)$.
\end{proposition}

\textbf{Proof.} 
By Schauder estimates, it suffices to prove that 
\[\sup \, (\varepsilon + |y|)^{1+\nu} \, |a(x,y)| \leq C \, \sup \, (\varepsilon + |y|)^{3+\nu} 
\, |\mathbb{L}_B a(x,y)|.\] 
Suppose that this estimate fails. Then there exists a sequence of $1$-forms $a^{(j)} \in 
\mathcal{E}_{1+\nu}^{2,\gamma}(\mathbb{R}^{n-4} \times \mathbb{R}^4)$ such that 
\[\sup \, (\varepsilon + |y|)^{1+\nu} \, |a^{(j)}(x,y)| = 1\] 
and 
\[\sup \, (\varepsilon + |y|)^{3+\nu} \, |\mathbb{L}_B a^{(j)}(x,y)| \to 0.\] 
Then there exists a sequence of points $(x_j,y_j) \in \mathbb{R}^{n-4} \times \mathbb{R}^4$ such 
that 
\[\sup \, (\varepsilon + |y_j|)^{1+\nu} \, |a^{(j)}(x_j,y_j)| \geq \frac{1}{2}.\] 
There are two possibilities: \\ 

(i) Suppose that the sequence $|y_j|$ is bounded. After passing to a subsequence, we may assume 
that the sequence $a^{(j)}$ converges to a $1$-form $a \in \mathcal{E}_{1+\nu}^{2,\gamma}
(\mathbb{R}^{n-4} \times \mathbb{R}^4)$ such that 
\[\sup \, (\varepsilon + |y|)^{1+\nu} \, |a(x,y)| \leq 1\] 
and 
\[\mathbb{L}_B a = 0.\] 
Using Proposition 3.3, we conclude that $a = 0$. This is a contradiction. \\

(ii) We now assume that $|y_j| \to \infty$. Let 
\[\tilde{a}^{(j)}(x,y) = |y_j|^{1+\nu} \, a^{(j)}(x + x_j,|y_j| \, y).\] 
After passing to a subsequence, we may assume that the sequence $\tilde{a}^{(j)}$ converges to a 
$1$-form $\tilde{a}$ such that 
\[\sup \, |y|^{1+\nu} \, |\tilde{a}(x,y)| \leq 1\] 
and 
\[d^* d \tilde{a} + dd^* \tilde{a} = 0.\] 
Thus, we conclude that $\tilde{a} = 0$. This is a contradiction. \\

\begin{proposition}
Let $0 < \nu < 1$. Assume that $b \in \mathcal{C}_{3+\nu}^\gamma(\mathbb{R}^4)$ satisfies 
\[\int_{\mathbb{R}^4} \langle b,F_B(X,\cdot) \rangle = 0\] 
for all $x \in \mathbb{R}^{n-4}$ and all vector fields of the form 
\[X = \varepsilon \, w_\rho \, e_\rho^\perp + \mu \, y_\rho \, e_\rho^\perp + r_{\rho\sigma} \, 
y_\sigma \, e_\rho^\perp\] 
with $w \in \mathbb{R}^4$, $\mu \in \mathbb{R}$, and $r \in \Lambda_+^2 \mathbb{R}^4$. Moreover, 
let $\eta \in \mathcal{S}(\mathbb{R}^{n-4})$. Then there exists a $1$-form $a \in \mathcal{E}
_{1+\nu}^{2,\gamma}(\mathbb{R}^{n-4} \times \mathbb{R}^4)$ such that 
\[\mathbb{L}_B a = \eta(x) \, b(y).\]
\end{proposition}

\textbf{Proof.} 
Let 
\[\eta(x) = \int_{\mathbb{R}^{n-4}} e^{ix\xi} \, \hat{\eta}(\xi) \, d\xi.\] 
For every $\xi \in \mathbb{R}^{n-4}$, there exists a $1$-form $\hat{a}(\xi,\cdot) \in \mathcal
{C}_{1+\nu}^{2,\gamma}(\mathbb{R}^4)$ such that 
\[\sum_{\beta=1}^4 D_{B,\beta} D_{B,\beta} \hat{a}_\alpha(\xi,y) + 2 \sum_{\beta=1}^4 [F_{B,
\alpha\beta},\hat{a}_\beta(\xi,y)] - |\xi|^2 \, \hat{a}_\alpha(\xi,y) = -b_\alpha(y)\] 
for $1 \leq \alpha \leq 4$ and 
\[\sum_{\beta=1}^4 D_{B,\beta} D_{B,\beta} \hat{a}_i(\xi,y) - |\xi|^2 \, \hat{a}_i(\xi,y) = -b_i
(y)\] 
for $1 \leq i \leq n - 4$. Furthermore, $\hat{a}(\xi,\cdot)$ satisfies 
\[\int_{\mathbb{R}^4} \langle \hat{a}(\xi,\cdot),F_B(X,\cdot) \rangle = 0\] 
for all $x \in \mathbb{R}^{n-4}$ and all vector fields of the form 
\[X = \varepsilon \, w_\rho \, e_\rho^\perp + \mu \, y_\rho \, e_\rho^\perp + r_{\rho\sigma} \, 
y_\sigma \, e_\rho^\perp = 0\] 
with $w \in \mathbb{R}^4$, $\mu \in \mathbb{R}$, and $r \in \Lambda_+^2 \mathbb{R}^4$. We now 
define a $1$-form $a \in \mathcal{E}_{1+\nu}^{2,\gamma}(\mathbb{R}^{n-4} \times \mathbb{R}^4)$ 
by 
\[a_\alpha(x,y) = \int_{\mathbb{R}^{n-4}} e^{ix\xi} \, \hat{\eta}(\xi) \, \hat{a}_\alpha(\xi,y) 
\, d\xi\] 
for $1 \leq \alpha \leq 4$ and 
\[a_i(x,y) = \int_{\mathbb{R}^{n-4}} e^{ix\xi} \, \hat{\eta}(\xi) \, \hat{a}_i(\xi,y) \, d\xi\] 
for $1 \leq i \leq n-4$. Then the $1$-form $a$ satisfies 
\[\sum_{j=1}^{n-4} \partial_j \partial_j a_\alpha(x,y) + \sum_{\beta=1}^4 D_{B,\beta} D_{B,
\beta} a_\alpha + 2 \sum_{\beta=1}^4 [F_{B,\alpha\beta},a_\beta] = -\eta(x) \, b_\alpha(y)\] 
for $1 \leq \alpha \leq 4$ and 
\[\sum_{j=1}^{n-4} \partial_j \partial_j a_i(x,y) + \sum_{\beta=1}^4 D_{B,\beta} D_{B,\beta} a_i 
= -\eta(x) \, b_i(y)\] 
for $1 \leq i \leq n - 4$. Thus, we conclude that $a \in \mathcal{E}_{1+\nu}^{2,\gamma}(\mathbb
{R}^{n-4} \times \mathbb{R}^4)$ and $\mathbb{L}_B a = \eta(x) \, b(y)$. \\

\begin{corollary}
Let $0 < \nu < 1$, and suppose that $b \in \mathcal{E}_{3+\nu}^\gamma(\mathbb{R}^{n-4} \times 
\mathbb{R}^4)$ has compact support. Then there exists a $1$-form $a \in \mathcal{E}_{1+\nu}
^{2,\gamma}(\mathbb{R}^{n-4} \times \mathbb{R}^4)$ such that 
\[\|a\|_{\mathcal{C}_{1+\nu}^{2,\gamma}(\mathbb{R}^{n-4} \times \mathbb{R}^4)} \leq C \, \|b\|
_{\mathcal{C}_{3+\nu}^\gamma(\mathbb{R}^{n-4} \times \mathbb{R}^4)}\] 
and 
\[\mathbb{L}_B a = b.\] 
\end{corollary}

\textbf{Proof.} 
It follows from Proposition 3.4 that the range of the operator $\mathbb{L}_B: \mathcal{E}_{1+
\nu}^{2,\gamma}(\mathbb{R}^{n-4} \times \mathbb{R}^4) \to \mathcal{E}_{3+\nu}^\gamma(\mathbb{R}
^{n-4} \times \mathbb{R}^4)$ is a closed subspace of the Banach space $\mathcal{E}_{3+\nu}
^\gamma(\mathbb{R}^{n-4} \times \mathbb{R}^4)$. By Proposition 3.5, it contains all $1$-forms of 
the form $\eta(x) \, b(y)$, where $\eta \in \mathcal{S}(\mathbb{R}^{n-4})$ and $b \in \mathcal
{C}_{3+\nu}^\gamma(\mathbb{R}^4)$ satisfies 
\[\int_{\mathbb{R}^4} \langle b,F_B(X,\cdot) \rangle = 0\] 
for all $x \in \mathbb{R}^{n-4}$ and all vector fields of the form 
\[X = \varepsilon \, w_\rho \, e_\rho^\perp  + \mu \, y_\rho \, e_\rho^\perp + r_{\rho\sigma} \, 
y_\sigma \, e_\rho^\perp\] 
with $w \in \mathbb{R}^4$, $\mu \in \mathbb{R}$, and $r \in \Lambda_+^2 \mathbb{R}^4$. The 
assertion follows now by approximation. \\

\begin{proposition}
Let $0 < \nu < 1$. Suppose that $b \in \mathcal{E}_{3+\nu}^\gamma(\mathbb{R}^{n-4} \times 
\mathbb{R}^4)$ is supported in the set $\{(x,y) \in \mathbb{R}^{n-4} \times \mathbb{R}^4: |x| 
\leq \delta, \, |y| \leq 2\delta^4\}$. Then there exists a $1$-form $a \in \mathcal{C}_{1+\nu}
^{2,\gamma}(\mathbb{R}^{n-4} \times \mathbb{R}^4)$ such that $a$ is supported in $\{(x,y) \in 
\mathbb{R}^{n-4} \times \mathbb{R}^4: |x| \leq 2\delta, \, |y| \leq 2\delta^2\}$, 
\[\|a\|_{\mathcal{C}_{1+\nu}^{2,\gamma}(\mathbb{R}^{n-4} \times \mathbb{R}^4)} \leq C \, \|b\|
_{\mathcal{C}_{3+\nu}^\gamma(\mathbb{R}^{n-4} \times \mathbb{R}^4)}\] 
and 
\[\|\mathbb{L}_B a - b\|_{\mathcal{C}_{3+\nu}^\gamma(\{(x,y) \in \mathbb{R}^{n-4} \times \mathbb
{R}^4: |y| \leq 2\delta^4\})} \leq C \, \delta \, \|b\|_{\mathcal{C}_{3+\nu}^\gamma(\mathbb{R}
^{n-4} \times \mathbb{R}^4)}\] 
and 
\[\|\mathbb{L}_B a - b\|_{\mathcal{C}_{3+\nu}^\gamma(\mathbb{R}^{n-4} \times \mathbb{R}^4)} \leq 
C \, |\log \delta|^{-1} \, \|b\|_{\mathcal{C}_{3+\nu}^\gamma(\mathbb{R}^{n-4} \times \mathbb{R}
^4)}.\] 
\end{proposition}

\textbf{Proof.} 
By Corollary 3.6, there exists a $1$-form $a \in \mathcal{E}_{1+\nu}^{2,\gamma}(\mathbb{R}^{n-4} 
\times \mathbb{R}^4)$ such that 
\[\|a\|_{\mathcal{C}_{1+\nu}^{2,\gamma}(\mathbb{R}^{n-4} \times \mathbb{R}^4)} \leq C \, \|b\|
_{\mathcal{C}_{3+\nu}^\gamma(\mathbb{R}^{n-4} \times \mathbb{R}^4)}\] 
and 
\[\mathbb{L}_B a = b.\] 
Let $\zeta$ be a cut-off function on $\mathbb{R}^{n-4}$ such that $\zeta(x) = 1$ for $|x| \leq 
\delta$, $\zeta(x) = 0$ for $|x| \geq 2\delta$, and 
\[\sup \, \delta \, |\nabla \zeta| + \sup \, \delta^2 \, |\nabla^2 \zeta| \leq C.\] 
Furthermore, let $\eta$ be a cut-off function on $\mathbb{R}^{n-4}$ satisfying $\eta(y) = 1$ for 
$|y| \leq 2\delta^4$, $\eta(y) = 0$ for $|y| \geq 2\delta^2$, and 
\[\sup \, |y| \, |\nabla \eta| + \sup \, |y|^2 \, |\nabla^2 \eta| \leq C \, |\log \delta|^{-1}.
\] 
Then we have the estimates 
\[\|\eta \, \zeta \, a\|_{\mathcal{C}_{1+\nu}^{2,\gamma}(\mathbb{R}^{n-4} \times \mathbb{R}^4)} 
\leq C \, \|b\|_{\mathcal{C}_{3+\nu}^\gamma(\mathbb{R}^{n-4} \times \mathbb{R}^4)}\] 
and 
\begin{align*} 
&\|\mathbb{L}_B(\zeta \, a) - b\|_{\mathcal{C}_{3+\nu}^\gamma(\{(x,y) \in \mathbb{R}^{n-4} 
\times \mathbb{R}^4: |y| \leq 2\delta^4\})} \\ 
&= \|\mathbb{L}_B (\zeta \, a) - \zeta \, \mathbb{L}_B a\|_{\mathcal{C}_{3+\nu}^\gamma(\{(x,y) 
\in \mathbb{R}^{n-4} \times \mathbb{R}^4: |y| \leq 2\delta^4\})} \\ 
&\leq C \, \delta \, \|a\|_{\mathcal{C}_{1+\nu}^{1,\gamma}(\mathbb{R}^{n-4} \times \mathbb{R}
^4)} \\ 
&\leq C \, \delta \, \|b\|_{\mathcal{C}_{3+\nu}^\gamma(\mathbb{R}^{n-4} \times \mathbb{R}^4)} 
\end{align*} 
and 
\begin{align*} 
&\|\mathbb{L}_B(\eta \, \zeta \, a) - b\|_{\mathcal{C}_{3+\nu}^\gamma(\mathbb{R}^{n-4} \times 
\mathbb{R}^4)} \\ 
&= \|\mathbb{L}_B (\eta \, \zeta \, a) - \eta \, \zeta \, \mathbb{L}_B a\|_{\mathcal{C}_{3+\nu}
^\gamma(\mathbb{R}^{n-4} \times \mathbb{R}^4)} \\ 
&\leq C \, |\log \delta|^{-1} \, \|a\|_{\mathcal{C}_{1+\nu}^{1,\gamma}(\mathbb{R}^{n-4} \times 
\mathbb{R}^4)} \\ 
&\leq C \, |\log \delta|^{-1} \, \|b\|_{\mathcal{C}_{3+\nu}^\gamma(\mathbb{R}^{n-4} \times 
\mathbb{R}^4)}. 
\end{align*} 
From this the assertion follows. \\

\section{Construction of the approximate solutions}

In this section, we describe the construction of the approximate solutions. To this end, we 
assume that the normal bundle $NS$ can be endowed with a $SU(2)$-structure $(J,\omega)$. Here, 
$J$ is a complex structure and $\omega$ is a complex volume form on $NS$. \\

Let $\nabla' = \nabla + \theta$ be a connection on the normal bundle $NS$ such that $\theta$ is 
a $1$-form with values in the Lie algebra $\Lambda_+^2 NS$ and $(J,\omega)$ is parallel with 
respect to the connection $\nabla'$. The $1$-form $\theta$ is uniquely determined by the 
covariant derivative of the pair $(J,\omega)$ with respect to the Levi-Civita connection 
$\nabla$. Since $(J,\omega)$ is parallel with respect to $\nabla'$, the connection induced by 
$\nabla'$ on the bundle $\Lambda_+^2 NS$ is flat. \\

The connection $\nabla'$ induces a splitting of the tangent space $TNS$ into horizontal and 
vertical subspaces. Let $\{e_i': 1 \leq i \leq n-4\}$ be an orthonormal basis for the horizontal 
subspace with respect to $\nabla'$, and let $\{e_\alpha^\perp: 1 \leq \alpha \leq 4\}$ be a 
$SU(2)$ basis for the vertical subspace $V$. \\

In the first step, we define a connection on the pull-back bundle $\pi^* NS$ of the normal 
bundle under the natural projection $\pi: NS \to S$. Since we may identify a neighborhood of $S$ 
in $M$ with a neighborhood of the zero section in $NS$, this gives a connection on a small 
neighborhood of $S$ in $M$. In the second step, we show that this connection can be extended to 
the whole of $M$ using suitable cut-off functions. \\

The glueing data consist of a set $(v,\lambda,J,\omega)$, where $v$ is a section of the normal 
bundle $NS$, $\lambda$ is a positive function on $S$, and $(J,\omega)$ is a $SU(2)$ structure on 
the normal bundle $NS$. As in Section 3, let $\{\mathfrak{i},\mathfrak{j},\mathfrak{k}\}$ be a 
basis for the Lie algebra $\mathfrak{su}(NS)$ such that  
\begin{align*} 
&\mathfrak{i}(e_1^\perp) = -e_2^\perp, \quad \mathfrak{i}(e_2^\perp) = e_1^\perp, \quad 
\mathfrak{i}(e_3^\perp) = e_4^\perp, \quad \mathfrak{i}(e_4^\perp) = -e_3^\perp, \\ 
&\mathfrak{j}(e_1^\perp) = -e_3^\perp, \quad \mathfrak{j}(e_2^\perp) = -e_4^\perp, \quad 
\mathfrak{j}(e_3^\perp) = e_1^\perp, \quad \mathfrak{j}(e_4^\perp) = e_2^\perp, \\ 
&\mathfrak{k}(e_1^\perp) = -e_4^\perp, \quad \mathfrak{k}(e_2^\perp) = e_3^\perp, \quad 
\mathfrak{k}(e_3^\perp) = -e_2^\perp, \quad \mathfrak{k}(e_4^\perp) = e_1^\perp. 
\end{align*} 
We consider a connection of the form $D_A = \nabla' + A$. The vertical components of $A$ are 
defined by 
\begin{align*} 
A(e_1^\perp) &= \frac{-(y - \varepsilon v)_2 \, \mathfrak{i} - (y - \varepsilon v)_3 \, 
\mathfrak{j} - (y - \varepsilon v)_4 \, \mathfrak{k}}{\varepsilon^2 \lambda^2 + |y - \varepsilon 
v|^2} \\ 
A(e_2^\perp) &= \frac{(y - \varepsilon v)_1 \, \mathfrak{i} - (y - \varepsilon v)_4 \, \mathfrak
{j} + (y - \varepsilon v)_3 \, \mathfrak{k}}{\varepsilon^2 \lambda^2 + |y - \varepsilon v|^2} \\ 
A(e_3^\perp) &= \frac{(y - \varepsilon v)_4 \, \mathfrak{i} + (y - \varepsilon v)_1 \, \mathfrak
{j} - (y - \varepsilon v)_2 \, \mathfrak{k}}{\varepsilon^2 \lambda^2 + |y - \varepsilon v|^2} \\ 
A(e_4^\perp) &= \frac{-(y - \varepsilon v)_3 \, \mathfrak{i} + (y - \varepsilon v)_2 \, 
\mathfrak{j} + (y - \varepsilon v)_1 \, \mathfrak{k}}{\varepsilon^2 \lambda^2 + |y - \varepsilon 
v|^2}.
\end{align*} 
Since the basic instanton on $\mathbb{R}^4$ is $SU(2)$-equivariant, this definition is 
independent of the choice of $SU(2)$-frame $\{e_\alpha^\perp: 1 \leq \alpha \leq 4\}$. 
Furthermore, the horizontal components of $A$ are defined by 
\[A(e_i') = -\varepsilon \, \nabla_i' v_\rho \, A(e_\rho^\perp) - \lambda^{-1} \, \nabla_i 
\lambda \, (y - \varepsilon v)_\rho \, A(e_\rho^\perp)\] 
for $1 \leq i \leq n-4$. \\

\begin{lemma}
The curvature of $A$ is given by 
\[F_A(e_i',e_\alpha^\perp) = -\Big ( \varepsilon \, \nabla_i' v_\rho + \lambda^{-1} \, \nabla_i 
\lambda \, (y - \varepsilon v)_\rho \Big ) \, F_A(e_\rho^\perp,e_\alpha^\perp)\] 
and 
\begin{align*} 
F_A(e_i',e_j') 
&= \Big ( \varepsilon \, \nabla_i' v_\rho + \lambda^{-1} \, \nabla_i \lambda \, (y - \varepsilon 
v)_\rho \Big ) \\ 
&\cdot \Big ( \varepsilon \, \nabla_j' v_\sigma + \lambda^{-1} \, \nabla_j \lambda \, (y - 
\varepsilon v)_\sigma \Big ) \, F_A(e_\rho^\perp,e_\sigma^\perp) \\ 
&+ C_{ij} + A \big ( C_{ij} \, (y - \varepsilon v) \big ), 
\end{align*} 
where $C_{ij} \in \Lambda_-^2 NS$ is the curvature of the connection $\nabla'$.
\end{lemma}

\textbf{Proof.} 
By definition of $A$, we have 
\begin{align*} 
F_A(e_i',e_\alpha^\perp) 
&= \nabla_{e_i'}' A(e_\alpha^\perp) - \nabla_{e_\alpha^\perp} A(e_i') + [A(e_i'),A(e_\alpha
^\perp)] \\ 
&= -\varepsilon \, \nabla_i' v_\rho \, \Big ( \nabla_{e_\rho^\perp} A(e_\alpha^\perp) - \nabla
_{e_\alpha^\perp} A(e_\rho^\perp) + [A(e_\rho^\perp),A(e_\alpha^\perp)] \Big ) \\ 
&- \lambda^{-1} \, \nabla_i \lambda \, (y - \varepsilon v)_\rho \, \Big ( \nabla_{e_\rho^\perp} 
A(e_\alpha^\perp) - \nabla_{e_\alpha^\perp} A(e_\rho^\perp) + [A(e_\rho^\perp),A(e_\alpha
^\perp)] \Big ) \\ 
&= -\varepsilon \, \nabla_i' v_\rho \, F_A(e_\rho^\perp,e_\alpha^\perp) \\ 
&- \lambda^{-1} \, \nabla_i \lambda \, (y - \varepsilon v)_\rho \, F_A(e_\rho^\perp,e_\alpha
^\perp) 
\end{align*} 
and 
\begin{align*} 
F_A(e_i',e_j') 
&= C_{ij} + \nabla_{e_i'}' A(e_j') - \nabla_{e_j'}' A(e_i') + [A(e_i'),A(e_j')] - A([e_i',e_j']) 
\\ 
&= \nabla_{e_i'}' A(e_j') - \nabla_{e_j'}' A(e_i') + [A(e_i'),A(e_j')] + C_{ij} + A \big ( 
C_{ij} \, y \big ) \\ 
&= \varepsilon \, \nabla_i' v_\rho \, \varepsilon \, \nabla_j' v_\sigma \, \Big ( \nabla_{e_\rho
^\perp} A(e_\sigma^\perp) - \nabla_{e_\sigma^\perp} A(e_\rho^\perp) + [A(e_\rho^\perp),A(e
_\sigma^\perp)] \Big ) \\ 
&+ \varepsilon \, \nabla_i' v_\rho \, \lambda^{-1} \, \nabla_j \lambda \, (y - \varepsilon v)
_\sigma \, \Big ( \nabla_{e_\rho^\perp} A(e_\sigma^\perp) - \nabla_{e_\sigma^\perp} A(e_\rho
^\perp) + [A(e_\rho^\perp),A(e_\sigma^\perp)] \Big ) \\ 
&+ \lambda^{-1} \, \nabla_i \lambda \, (y - \varepsilon v)_\rho \, \varepsilon \, \nabla_j' v
_\sigma \, \Big ( \nabla_{e_\rho^\perp} A(e_\sigma^\perp) - \nabla_{e_\sigma^\perp} A(e_\rho
^\perp) + [A(e_\rho^\perp),A(e_\sigma^\perp)] \Big ) \\ 
&+ C_{ij} + A \big ( C_{ij} \, (y - \varepsilon v) \big ) \\ 
&= \varepsilon \, \nabla_i' v_\rho \, \varepsilon \, \nabla_j' v_\sigma \, F_A(e_\rho^\perp,e
_\sigma^\perp) \\ 
&+ \varepsilon \, \nabla_i' v_\rho \, \lambda^{-1} \, \nabla_j \lambda \, (y - \varepsilon v)
_\sigma \, F_A(e_\rho^\perp,e_\sigma^\perp) \\ 
&+ \lambda^{-1} \, \nabla_i \lambda \, (y - \varepsilon v)_\rho \, \varepsilon \, \nabla_j' v
_\sigma \, F_A(e_\rho^\perp,e_\sigma^\perp) \\ 
&+ C_{ij} + A \big ( C_{ij} \, (y - \varepsilon v) \big ). 
\end{align*}
This proves the assertion. \\

Let $\{e_i: 1 \leq i \leq n-4\}$ be an orthonormal basis for the horizontal subspace with 
respect to the Levi-Civita connection $\nabla$. Then we have the following result: \\

\begin{lemma}
The curvature of $A$ satisfies 
\[F_A(e_i,e_\alpha^\perp) = -\Big ( \varepsilon \, \nabla_i v_\rho + \lambda^{-1} \, \nabla_i 
\lambda \, (y - \varepsilon v)_\rho + \theta_{i,\rho\sigma} \, (y - \varepsilon v)_\sigma \Big ) 
\, F_A(e_\rho^\perp,e_\alpha^\perp)\] 
and 
\begin{align*} 
F_A(e_i,e_j) 
&= \Big ( \varepsilon \, \nabla_i v_\rho + \lambda^{-1} \, \nabla_i \lambda \, (y - \varepsilon 
v)_\rho + \theta_{i,\rho\alpha} \, (y - \varepsilon v)_\alpha \Big ) \\ 
&\cdot \Big ( \varepsilon \, \nabla_j v_\sigma + \lambda^{-1} \, \nabla_j \lambda \, (y - 
\varepsilon v)_\sigma + \theta_{j,\sigma\beta} \, (y - \varepsilon v)_\beta \Big ) \, F_A(e_\rho
^\perp,e_\sigma^\perp) \\ 
&+ C_{ij} + A \big ( C_{ij} \, (y - \varepsilon v) \big ), 
\end{align*} 
where $C_{ij} \in \Lambda_-^2 NS$ is the curvature of $\nabla'$.
\end{lemma}

\textbf{Proof.} 
Since $\nabla' = \nabla + \theta$, the orthonormal basis $\{e_i: 1 \leq i \leq n-4\}$ is related 
to the orthonormal basis $\{e_i': 1 \leq i \leq n-4\}$ by 
\[e_i' = e_i + \theta_{i,\rho\sigma} \, y_\sigma \, e_\rho^\perp.\] 
The assertion follows now from Lemma 4.1. \\

Using Lemma 2.1, we obtain the following result. \\

\begin{lemma}
Suppose that $\mu$ is constant and $r$ is a section of the vector bundle $\Lambda_+^2 NS$ such 
that $\nabla' r = 0$. Let 
\[u = (\varepsilon^2 \lambda^2  + |y - \varepsilon v|^2)^{-\frac{1}{2}} \, \big ( \mu \, 
(y - \varepsilon v)_\rho + r_{\rho\sigma} \, (y - \varepsilon v)_\sigma \big ) \, e_\rho
^\perp.\] 
Then the covariant derivative of $u$ satisfies the estimate 
\[\|D_A u\|_{\mathcal{C}_3^{1,\gamma}(M)} \leq C \, \varepsilon^2.\] 
\end{lemma}

\vspace{2mm}

Hence, as we move away from the submanifold $S$, the connection $A$ approaches a flat 
connection. Therefore, we can extend $A$ trivially to $M$. \\

Our aim is to derive estimates for $D_A^* F_A$ in $\mathcal{C}_3^\gamma(M)$. To this end, we 
assume that the glueing data $(v,\lambda,J,\omega)$ satisfy the estimates 
\[\|v\|_{\mathcal{C}^{2,\gamma}(S)} \leq K,\] 
\[\|\lambda\|_{\mathcal{C}^{2,\gamma}(S)} \leq K, \qquad \inf \lambda \geq 1,\] 
\[\|(J,\omega)\|_{\mathcal{C}^{2,\gamma}(S)} \leq K\] 
for some $K > 0$. All implicite constants will depend on $K$. \\

\begin{proposition}
If the set $(v,\lambda,J,\omega)$ is admissible, then we have the estimate 
\[\|D_A^* F_A\|_{\mathcal{C}_3^\gamma(M)} \leq C \, \varepsilon^2.\] 
\end{proposition}

\textbf{Proof.} 
Since $B$ is a Yang-Mills connection on $\mathbb{R}^4$, we have 
\[\sum_{\beta=1}^4 D_{A,e_\beta^\perp} F_A(e_\alpha^\perp,e_\beta^\perp) = 0.\] 
Using Proposition 2.4, we obtain 
\[\sum_{\beta=1}^4 D_{A,e_\beta^\perp} \big ( \nabla_i v_\rho \, F_A(e_\rho^\perp,e_\beta^\perp) 
\big ) = 0,\] 
\[\sum_{\beta=1}^4 D_{A,e_\beta^\perp} \big ( \lambda^{-1} \, \nabla_i \lambda \, (y - 
\varepsilon v)_\rho \, F_A(e_\rho^\perp,e_\beta^\perp) \big ) = 0,\] 
\[\sum_{\beta=1}^4 D_{A,e_\beta^\perp} \big ( \theta_{i,\rho\sigma} \, (y - \varepsilon v)
_\sigma \, F_A(e_\rho^\perp,e_\beta^\perp) \big ) = 0.\] 
From this it follows that 
\[\sum_{\beta=1}^4 D_{A,e_\beta^\perp} F_A(e_i,e_\beta^\perp) = 0.\] 
Therefore, we obtain 
\[\|D_A^{*_{g_0}} F_A\|_{\mathcal{C}_3^\gamma(M)} \leq C \, \varepsilon^2.\] 
Here, $g_0$ denotes the product metric on $NS$, i.e. 
\begin{align*} 
&g_0(e_i,e_j) = \delta_{ij} \\ 
&g_0(e_i,e_\alpha^\perp) = 0 \\ 
&g_0(e_\alpha^\perp,e_\beta^\perp) = \delta_{\alpha\beta}. 
\end{align*} 
Let $g$ be the pull-back of the Riemannian metric on $M$ under the exponential map $\exp: NS \to 
M$. Then the metric $g$ satisfies an asymptotic expansion of the form 
\begin{align*} 
&g(e_i,e_j) = \delta_{ij} + 2 \sum_{\rho=1}^4 h_{ij,\rho} \, y_\rho + O(|y|^2) \\ 
&g(e_i,e_\alpha^\perp) = O(|y|^2) \\ 
&g(e_\alpha^\perp,e_\beta^\perp) = \delta_{\alpha\beta} + O(|y|^2), 
\end{align*} 
where $h$ denotes the second fundamental form of $S$. In particular, the volume form of $g$ is 
related to the volume form of $g_0$ by 
\[\bigg ( \frac{\det g}{\det g_0} \bigg )^{\frac{1}{2}} = 1 + H_\rho \, y_\rho + O(|y|^2),\] 
where $H$ is the mean curvature vector of $S$. Since the mean curvature of $S$ is $0$, we obtain 
\[\bigg ( \frac{\det g}{\det g_0} \bigg )^{\frac{1}{2}} = 1 + O(|y|^2).\] 
Thus, we conclude that 
\[\|D_A^* F_A\|_{\mathcal{C}_3^\gamma(M)} \leq C \, \varepsilon^2\] 

\vspace{2mm}

\section{Estimates for the linearized operator in weighted H\"older spaces}

Our aim in this section is to analyze the mapping properties of the linearized operator $\mathbb
{L}_A: \Omega^1(M) \to \Omega^1(M)$. \\

\begin{proposition}
Suppose that $b \in \mathcal{C}_{3+\nu}^\gamma(M)$ is supported in the set $\{p \in M: \dist(p,
S) \leq 2\delta^4\}$ and satisfies 
\[\int_{NS_x} \sum_{\alpha=1}^4 \langle b(e_\alpha^\perp),F_A(X,e_\alpha^\perp) \rangle = 0\] 
for all $x \in S$ and all vector fields of the form 
\[X = \varepsilon \, w_\rho \, e_\rho^\perp + \mu \, (y - \varepsilon v)_\rho \, e_\rho^\perp 
+ r_{\rho\sigma} \, (y - \varepsilon v)_\sigma \, e_\rho^\perp\] 
with $w \in NS$, $\mu \in \mathbb{R}$, and $r \in \Lambda_+^2 NS$. Then there exists a $1$-form 
$a \in \mathcal{C}_{1+\nu}^{2,\gamma}(M)$ which is supported in the region $\{p \in M: \dist(p,
S) \leq 2\delta^2\}$ such that 
\[\|a\|_{\mathcal{C}_{1+\nu}^{2,\gamma}(M)} \leq C \, \|b\|_{\mathcal{C}_{3+\nu}^\gamma(M)}\] 
and 
\[\|\mathbb{L}_A a - b\|_{\mathcal{C}_{3+\nu}^\gamma(\{p \in M: \dist(p,S) \leq 2\delta^4\})} 
\leq C \, \delta \, \|b\|_{\mathcal{C}_{3+\nu}^\gamma(M)}\] 
and 
\[\|\mathbb{L}_A a - b\|_{\mathcal{C}_{3+\nu}^\gamma(M)} \leq C \, |\log \delta|^{-1} \, \|b\|
_{\mathcal{C}_{3+\nu}^\gamma(M)}.\]
\end{proposition}

\textbf{Proof.} 
Let $\{\zeta^{(j)}: 1 \leq j \leq j_0\}$ be a partition of unity on $S$ such that each function 
$\zeta^{(j)}$ is supported in a ball $B_\delta(p_j)$, and 
\[|\{1 \leq j \leq j_0: x \in B_{4\delta}(p_j)\}| \leq C\] 
for all $x \in S$ and some uniform constant $C$. For each $1 \leq j \leq j_0$, 
there exists a $1$-form $a^{(j)} \in \mathcal{C}_{1+\nu}
^{2,\gamma}(M)$ which is supported in the region $\{(x,y) \in NS: x \in B_{2\delta}(p_j), \, |y| 
\leq 2\delta^2\}$ such that 
\[\|a^{(j)}\|_{\mathcal{C}_{1+\nu}^{2,\gamma}(M)} \leq C \, \|\zeta^{(j)} \, b\|
_{\mathcal{C}_{3+\nu}^\gamma(M)}\] 
and 
\[\|\mathbb{L}_A a^{(j)} - \zeta^{(j)} \, b\|_{\mathcal{C}_{3+\nu}^\gamma(\{p \in M: \dist(p,S) 
\leq 2\delta^4\})} \leq C \, \delta \, \|\zeta^{(j)} \, b\|_{\mathcal{C}_{3+\nu}^\gamma(M)}\] 
and 
\[\|\mathbb{L}_A a^{(j)} - \zeta^{(j)} \, b\|_{\mathcal{C}_{3+\nu}^\gamma(M)} \leq C \, |\log 
\delta|^{-1} \, \|b\|_{\mathcal{C}_{3+\nu}^\gamma(M)}.\] 
We now define 
\[a = \sum_{j=1}^{j_0} a^{(j)}.\] 
Then we have the estimates 
\begin{align*} 
\|a\|_{\mathcal{C}_{1+\nu}^{2,\gamma}(M)} 
&\leq C \, \sup_{1 \leq j \leq j_0} \|a^{(j)}\|_{\mathcal{C}_{1+\nu}^{2,\gamma}(M)} \\ 
&\leq C \, \sup_{1 \leq j \leq j_0} \|\zeta^{(j)} \, b\|_{\mathcal{C}_{3+\nu}^\gamma(M)} \\ 
&\leq C \, \|b\|_{\mathcal{C}_{3+\nu}^\gamma(M)}, 
\end{align*} 
\begin{align*} 
\|\mathbb{L}_A a - b\|_{\mathcal{C}_{3+\nu}^\gamma(\{p \in M: \dist(p,S) \leq 2\delta^4\})} 
&\leq C \, \sup_{1 \leq j \leq j_0} \|\mathbb{L}_A a^{(j)} - \zeta^{(j)} \, b\|_{\mathcal{C}
_{3+\nu}^\gamma(\{p \in M: \dist(p,S) \leq 2\delta^4\})} \\ 
&\leq C \, \delta \, \sup_{1 \leq j \leq j_0} \|\zeta^{(j)} \, b\|_{\mathcal{C}_{3+\nu}^\gamma
(M)} \\ 
&\leq C \, \delta \, \|b\|_{\mathcal{C}_{3+\nu}^\gamma(M)}, 
\end{align*} 
\begin{align*} 
\|\mathbb{L}_A a - b\|_{\mathcal{C}_{3+\nu}^\gamma(M)} 
&\leq C \, \sup_{1 \leq j \leq j_0} \|\mathbb{L}_A a^{(j)} - \zeta^{(j)} \, b\|_{\mathcal{C}_{3+
\nu}^\gamma(M)} \\ 
&\leq C \, |\log \delta|^{-1} \, \sup_{1 \leq j \leq j_0} \|\zeta^{(j)} \, b\|_{\mathcal{C}_{3+
\nu}^\gamma(M)} \\ 
&\leq C \, |\log \delta|^{-1} \, \|b\|_{\mathcal{C}_{3+\nu}^\gamma(M)}. 
\end{align*} 
This proves the assertion. \\

\begin{proposition}
For every $b \in \mathcal{C}_{3+\nu}^\gamma(M)$, there exists a $1$-form $a \in \mathcal{C}_{1+
\nu}^{2,\gamma}(M)$ such that 
\[\|a\|_{\mathcal{C}_{1+\nu}^{2,\gamma}(M)} \leq C \, \|b\|_{\mathcal{C}_{3+\nu}^\gamma(M)}\] 
and 
\[d^* da + dd^* a = b.\]
\end{proposition}

\textbf{Proof.} 
Since $H^1(M) = 0$, the operator $d^* d + dd^*: \Omega^1(M) \to \Omega^1(M)$ is invertible. 
Hence, there exists a $1$-form $a$ such that 
\[d^* da + dd^* a = b.\] 
Therefore, it remains to show that 
\[\|a\|_{\mathcal{C}_{1+\nu}^{2,\gamma}(M)} \leq C \, \|d^* da + dd^* a\|_{\mathcal{C}_{3+\nu}
^\gamma(M)}.\] 
By Schauder estimates, it suffices to show that 
\[\sup \, (\varepsilon + \dist(p,S))^{1+\nu} \, |a| \leq C \, \sup \, (\varepsilon + \dist(p,S))
^{3+\nu} \, |d^* da + dd^* a|.\] 
Suppose that this estimate fails. Then there exists a sequence of positive real numbers 
$\varepsilon_j$ and a sequence of $1$-forms $a^{(j)} \in \mathcal{C}_{1+\nu}^{2,\gamma}(M)$ such 
that 
\[\sup \, (\varepsilon_j + \dist(p,S))^{1+\nu} \, |a^{(j)}| = 1\] 
and 
\[\sup \, (\varepsilon_j + \dist(p,S))^{3+\nu} \, |d^* da^{(j)} + dd^* a^{(j)}| \to 0.\] 
Then there exists a sequence of points $p_j \in M$ such that 
\[\sup \, (\varepsilon_j + \dist(p_j,S))^{1+\nu} \, |a^{(j)}(p_j)| \geq \frac{1}{2}.\] 
There are two possibilities: \\ 

(i) Suppose that $\dist(p_j,S)$ is bounded from below. After passing to a subsequence, we may 
assume that the sequence $a^{(j)}$ converges to a $1$-form $a \in \Omega^1(M)$ such that 
\[\sup \, \dist(p,S)^{1+\nu} \, |a| \leq 1\] 
and 
\[d^* da + dd^* a = 0.\] 
From this it follows that $a$ is smooth. Since the operator $d^* d + dd^*: \Omega^1(M) \to 
\Omega^1(M)$ has trivial kernel, it follows that $a = 0$. This is a contradiction. \\

(ii) We now assume that $\dist(p_j,S) \to 0$. After rescaling and taking the limit, we obtain a 
$1$-form $\tilde{a} \in \Omega^1(\mathbb{R}^{n-4} \times \mathbb{R}^4)$ such that 
\[\sup \, |y|^{1+\nu} \, |\tilde{a}| \leq 1\] 
and 
\[d^* d \tilde{a} + dd^* \tilde{a} = 0.\] 
Thus, we conclude that $\tilde{a} = 0$. This is a contradiction. \\

\begin{proposition}
Suppose that $b \in \mathcal{C}_{3+\nu}^\gamma(M)$ is supported in the region $\{p \in M: \dist
(p,S) \geq \delta^4\}$. Then there exists a $1$-form $a \in \mathcal{C}_{1+\nu}^{2,\gamma}(M)$ 
which is supported in the region $\{p \in M: \dist(p,S) \geq \delta^8\}$ such that 
\[\|a\|_{\mathcal{C}_{1+\nu}^{2,\gamma}(M)} \leq C \, \|b\|_{\mathcal{C}_{3+\nu}^\gamma(M)}\] 
and 
\[\|\mathbb{L}_A a - b\|_{\mathcal{C}_{3+\nu}^\gamma(M)} \leq C \, \Big ( |\log \delta|^{-1} + 
\delta^{-16} \, \varepsilon^2 \Big ) \, \|b\|_{\mathcal{C}_{3+\nu}^\gamma(M)}.\]
\end{proposition}

\textbf{Proof.} 
By Proposition 5.2, there exists a $1$-form $a \in \mathcal{C}_{1+\nu}^{2,\gamma}(M)$ such that 
\[\|a\|_{\mathcal{C}_{1+\nu}^{2,\gamma}(M)} \leq C \, \|b\|_{\mathcal{C}_{3+\nu}^\gamma(M)}\] 
and 
\[d^* da + dd^* a = b.\] 
Let now $\eta$ be a cut-off function such that $\eta(p) = 0$ for $\dist(p,S) \leq \delta^8$, 
$\eta(p) = 1$ for $\dist(p,S) \geq \delta^4$ and 
\[\sup \, \dist(p,S) \, |\nabla \eta| + \sup \, \dist(p,S)^2 \, |\nabla^2 \eta| \leq C \, |\log 
\delta|^{-1}.\] 
Then the $1$-form $\eta \, a$ is supported in the region $\{p \in M: \dist(p,S) \geq \delta^8\}$ 
and satisfies 
\begin{align*} 
&\|\mathbb{L}_A (\eta \, a) - b\|_{\mathcal{C}_{3+\nu}^\gamma(M)} \\ 
&= \|D_A^* D_A(\eta \, a) + D_A D_A^*(\eta \, a) + (-1)^n \, *[*F_A,\eta \, a] - \eta \, (d^* da 
+ dd^* a)\|_{\mathcal{C}_{3+\nu}^\gamma(M)} \\ 
&\leq \|D_A^* D_A(\eta \, a) + D_A D_A^*(\eta \, a) + (-1)^n \, *[*F_A,\eta \, a] - d^* d(\eta 
\, a) - dd^*(\eta \, a)\|_{\mathcal{C}_{3+\nu}^\gamma(M)} \\ 
&+ \|d^* d(\eta \, a) + dd^*(\eta \, a) - \eta \, (d^* da + dd^* a)\|_{\mathcal{C}_{3+\nu}
^\gamma(M)} \\ 
&\leq \|D_A^* D_A(\eta \, a) + D_A D_A^*(\eta \, a) + (-1)^n \, *[*F_A,\eta \, a] - d^* d(\eta 
\, a) - dd^*(\eta \, a)\|_{\mathcal{C}_{3+\nu}^\gamma(M)} \\ 
&+ \|d^* d(\eta \, a) + dd^*(\eta \, a) - \eta \, (d^* da + dd^* a)\|_{\mathcal{C}_{3+\nu}
^\gamma(M)} \\ 
&\leq C \, \delta^{-16} \, \varepsilon^2 \, \|a\|_{\mathcal{C}_{1+\nu}^{1,\gamma}(M)} + C \, 
|\log \delta|^{-1} \, \|a\|_{\mathcal{C}_{1+\nu}^{1,\gamma}(M)} \\ 
&\leq C \, \delta^{-16} \, \varepsilon^2 \, \|b\|_{\mathcal{C}_{3+\nu}^\gamma(M)} + C \, |\log 
\delta|^{-1} \, \|b\|_{\mathcal{C}_{3+\nu}^\gamma(M)}. 
\end{align*} 
This proves the assertion. \\

In the following, we will choose $\delta = \varepsilon^{\frac{1}{16}}$. Let $\kappa$ be a 
cut-off function such that $\kappa(p) = 1$ for $\dist(p,S) \leq \varepsilon^{\frac{1}{4}}$ and 
$\kappa(p) = 0$ for $\dist(p,S) \geq 2 \, \varepsilon^{\frac{1}{4}}$. \\

Let $\mathcal{E}_\nu^{k,\gamma}(M)$ be the set of all $b \in \Omega^1(M)$ such that $b \in 
\mathcal{C}_\nu^{k,\gamma}(M)$ and 
\[\int_{NS_x} \kappa \, \sum_{\alpha=1}^4 \langle b(e_\alpha^\perp),F_A(X,e_\alpha^\perp) 
\rangle = 0\] 
for all $x \in S$ and all vector fields of the form 
\[X = \varepsilon \, w_\rho \, e_\rho^\perp + \mu \, (y - \varepsilon v)_\rho \, e_\rho^\perp 
+ r_{\rho\sigma} \, (y - \varepsilon v)_\sigma \, e_\rho^\perp\] 
with $w \in NS_x$, $\mu \in \mathbb{R}$, and $r \in \Lambda_+^2 NS_x$. \\

We denote by $\Id - \mathbb{P}$ the fibrewise projection from $\mathcal{C}_\nu^\gamma(M)$ to the 
subspace $\mathcal{E}_\nu^\gamma(M)$. Hence, if $b$ is a $1$-form, then the projection $\mathbb
{P} b$ is of the form 
\[\mathbb{P} b(e_\alpha^\perp) = \kappa \, \big ( \varepsilon \, w_\rho + \mu \, (y - 
\varepsilon v)_\rho + r_{\rho\sigma} \, (y - \varepsilon v)_\sigma \big ) \, F_A(e_\rho^\perp,
e_\alpha^\perp)\] 
for some $w \in NS$, $\mu \in \mathbb{R}$, and $r \in \Lambda_+^2 NS$. Let $\Pi$ be the linear 
operator which assigns to every $1$-form $b$ the triplet 
\[\Pi b = (w,\mu,r) \in NS \oplus \mathbb{R} \oplus \Lambda_+^2 NS.\] 
We shall need the following estimate for the operator norm of the projection operator $\mathbb
{P}$. \\

\begin{proposition}
For every $1$-form $b \in \mathcal{C}_{3+\nu}^\gamma(M)$, we have the estimates 
\[\|\Pi b\|_{\mathcal{C}^\gamma(S)} \leq C \, \varepsilon^{-2-\nu-\gamma} \, \|b\|_{\mathcal{C}
_{3+\nu}^\gamma(M)}\] 
and 
\[\|\mathbb{P} b\|_{\mathcal{C}_{3+\nu}^\gamma(M)} \leq C \, \varepsilon^{-\nu-\gamma} \, \|b\|
_{\mathcal{C}_{3+\nu}^\gamma(M)}.\] 
\end{proposition}

\textbf{Proof.} 
Without loss of generality, we may assume that 
\[\|b\|_{\mathcal{C}_{3+\nu}^\gamma(M)} \leq 1.\] 
Consider a point $x \in S$ and let $X$ be a vector field of the form 
\[X = \varepsilon \, w_\rho \, e_\rho^\perp + \mu \, (y - \varepsilon v)_\rho \, e_\rho^\perp 
+ r_{\rho\sigma} \, (y - \varepsilon v)_\sigma \, e_\rho^\perp.\] 
Using the estimate 
\[\sup \, (\varepsilon + |y|)^{3+\nu} \, |b(x,y)| \leq 1,\] 
we obtain 
\[\bigg | \int_{NS_x} \kappa \, \sum_{\alpha=1}^4 \langle b(e_\alpha^\perp),F_A(X,e_\alpha
^\perp) \rangle \bigg | \leq C \, \varepsilon^{-\nu} \, (|w| + |\mu| + |r|).\] 
From this it follows that 
\[\sup \, |\Pi b(x)| \leq C \, \varepsilon^{-2-\nu}.\] 
This implies 
\[\sup_{4 |x_1 - x_2| \geq \varepsilon} \frac{|\Pi b(x_1) - \Pi b(x_2)|}{|x_1 - x_2|^\gamma} 
\leq C \, \varepsilon^{-2-\nu-\gamma}.\] 
Using the estimate 
\[\sup_{4 |x_1 - x_2| \leq \varepsilon} \, (\varepsilon + |y|)^{3+\nu} \frac{|b(x_1,y) - b(x_2,
y)|}{|x_1 - x_2|^\gamma} \leq \varepsilon^{-\gamma},\] 
we deduce that 
\[\sup \, \frac{|\Pi b(x_1) - \Pi b(x_2)|}{|x_1 - x_2|^\gamma} \leq C \, \varepsilon^{-2-\nu-
\gamma}.\] 
Thus, we conclude that 
\[\|\Pi b\|_{\mathcal{C}^\gamma(S)} \leq C \, \varepsilon^{-2-\nu-\gamma},\] 
hence 
\[\|\mathbb{P} b\|_{\mathcal{C}_{3+\nu}^\gamma(M)} \leq C \, \varepsilon^{-\nu-\gamma}.\] 
This proves the assertion. \\

\begin{proposition}
For every $b \in \mathcal{E}_{3+\nu}^\gamma(M)$ there exists a $1$-form $a \in \mathcal{C}_{1+
\nu}^{2,\gamma}(M)$ such that 
\[\|a\|_{\mathcal{C}_{1+\nu}^{2,\gamma}(M)} \leq C \, \|b\|_{\mathcal{C}_{3+\nu}^\gamma(M)}\] 
and 
\[\|\mathbb{L}_A a - b\|_{\mathcal{C}_{3+\nu}^\gamma(\{p \in M: \dist(p,S) \leq \varepsilon
^{\frac{1}{2}}\})} \leq C \, \varepsilon^{\frac{1}{16}} \, \|b\|_{\mathcal{C}_{3+\nu}^\gamma(M)}
\] 
and 
\[\|\mathbb{L}_A a - b\|_{\mathcal{C}_{3+\nu}^\gamma(M)} \leq C \, |\log \varepsilon|^{-1} \, 
\|b\|_{\mathcal{C}_{3+\nu}^\gamma(M)}.\]
\end{proposition}

\textbf{Proof.} 
Apply Proposition 5.1 to $\kappa \, b$ and Proposition 5.3 to $(1 - \kappa) \, b$. \\

\begin{proposition}
For every $b \in \mathcal{E}_{3+\nu}^\gamma(M)$ there exists a $1$-form $a \in \mathcal{C}_{1+
\nu}^{2,\gamma}(M)$ such that 
\[\|a\|_{\mathcal{C}_{1+\nu}^{2,\gamma}(M)} \leq C \, \|b\|_{\mathcal{C}_{3+\nu}^\gamma(M)}\] 
and 
\[(\Id - \mathbb{P}) \, \mathbb{L}_A a = b.\] 
Furthermore, $a$ satisfies the estimate 
\[\|\Pi \, \mathbb{L}_A a\|_{\mathcal{C}^\gamma(S)} \leq C \, \varepsilon^{-2+\frac{1}{32}} \, 
\|b\|_{\mathcal{C}_{3+\nu}^\gamma(M)}.\] 
\end{proposition}

\textbf{Proof.} 
By Proposition 5.5, there exists an operator $\mathbb{S}: \mathcal{E}_{3+\nu}^\gamma(M) \to 
\mathcal{C}_{1+\nu}^{2,\gamma}(M)$ such that 
\[\|\mathbb{S} b\|_{\mathcal{C}_{1+\nu}^{2,\gamma}(M)} \leq C \, \|b\|_{\mathcal{C}_{3+\nu}
^\gamma(M)}\] 
and 
\[\|\mathbb{L}_A \, \mathbb{S} b - b\|_{\mathcal{C}_{3+\nu}^\gamma(\{p \in M: \dist(p,S) \leq 
\varepsilon^{\frac{1}{2}}\})} \leq C \, \varepsilon^{\frac{1}{16}} \, \|b\|_{\mathcal{C}_{3+\nu}
^\gamma(M)}\] 
and 
\[\|\mathbb{L}_A \, \mathbb{S} b - b\|_{\mathcal{C}_{3+\nu}^\gamma(M)} \leq C \, |\log 
\varepsilon|^{-1} \, \|b\|_{\mathcal{C}_{3+\nu}^\gamma(M)}.\]
This implies 
\[\|\Pi \, \mathbb{L}_A \, \mathbb{S} b\|_{\mathcal{C}^\gamma(S)} = \|\Pi(\mathbb{L}_A \, 
\mathbb{S} b - b)\|_{\mathcal{C}^\gamma(S)} \leq C \, \varepsilon^{-2+\frac{1}{16}-\nu-\gamma} 
\, \|b\|_{\mathcal{C}_{3+\nu}^\gamma(M)}.\] 
From this it follows that 
\[\|(\Id - \mathbb{P}) \, \mathbb{L}_A \, \mathbb{S} b - b\|_{\mathcal{C}_{3+\nu}^\gamma(M)} 
\leq C \, |\log \varepsilon|^{-1} \, \|b\|_{\mathcal{C}_{3+\nu}^\gamma(M)}.\] 
Therefore, the operator $(\Id - \mathbb{P}) \, \mathbb{L}_A \, \mathbb{S}: \mathcal{E}_{3+\nu}
^\gamma(M) \to \mathcal{E}_{3+\nu}^\gamma(M)$ is invertible. Hence, if we define 
\[a = \mathbb{S} \, \big [ (\Id - \mathbb{P}) \, \mathbb{L}_A \, \mathbb{S} \big ]^{-1} \, b,\] 
then $a$ satisfies 
\[\|a\|_{\mathcal{C}_{1+\nu}^{2,\gamma}(M)} \leq C \, \|b\|_{\mathcal{C}_{3+\nu}^\gamma(M)}\] 
and 
\[(\Id - \mathbb{P}) \, \mathbb{L}_A a = b.\] 
This proves the assertion. \\

\section{The nonlinear problem}

\begin{proposition}
For every approximate solution $A$, there exists a nearby connection $\tilde{A} = A + a$ such 
that 
\[\|a\|_{C_{1+\nu}^{2,\gamma}(M)} \leq C \, \varepsilon^{2-\nu-\gamma}\] 
and 
\[(\Id - \mathbb{P}) \, (D_{\tilde{A}}^* F_{\tilde{A}} + D_{\tilde{A}} D_{\tilde{A}}^* a) = 0.\] 
Furthermore, $a$ satisfies the estimate 
\[\|\Pi \, \mathbb{L}_A a\|_{\mathcal{C}^\gamma(S)} \leq C \, \varepsilon^{\frac{1}{32}}.\] 
\end{proposition}

\textbf{Proof.} 
The connection $\tilde{A} = A + a$ satisfies 
\[D_{\tilde{A}}^* F_{\tilde{A}} + D_{\tilde{A}} D_{\tilde{A}}^* a = D_A^* F_A + L_A a + D_A D_A
^* a + Q(a),\] 
where $Q(a)$ contains only quadratic and cubic terms. This implies 
\[D_{\tilde{A}}^* F_{\tilde{A}} + D_{\tilde{A}} D_{\tilde{A}}^* a = D_A^* F_A + \mathbb{L}_A a + 
Q(a).\] 
According to Proposition 5.6, there exists an operator $\mathbb{G}: \mathcal{E}_{3+\nu}^\gamma
(M) \to \mathcal{C}_{1+\nu}^{2,\gamma}(M)$ such that 
\[\|\mathbb{G} b\|_{\mathcal{C}_{1+\nu}^{2,\gamma}(M)} \leq C \, \|b\|_{\mathcal{C}_{3+\nu}
^\gamma(M)}\] 
and 
\[(\Id - \mathbb{P}) \, \mathbb{L}_A \, \mathbb{G} = \Id.\] 
We now define a mapping $\Phi: \mathcal{C}_{1+\nu}^{2,\gamma}(M) \to \mathcal{C}_{1+\nu}^{2,
\gamma}(M)$ by 
\[\Phi(a) = -\mathbb{G} \, (\Id - \mathbb{P}) \, (D_A^* F_A) - \mathbb{G} \, (\Id - \mathbb{P}) 
\, Q(a).\] 
Then we have the estimate 
\begin{align*} 
\|\Phi(a)\|_{\mathcal{C}_{1+\nu}^{2,\gamma}(M)} 
&\leq C \, \|(\Id - \mathbb{P}) \, (D_A^* F_A)\|_{\mathcal{C}_{3+\nu}^\gamma(M)} + C \, \| (\Id 
- \mathbb{P}) \, Q(a)\|_{\mathcal{C}_{3+\nu}^\gamma(M)} \\ 
&\leq C \, \varepsilon^{-\nu-\gamma} \, \|D_A^* F_A\|_{\mathcal{C}_{3+\nu}^\gamma(M)} + C \, 
\varepsilon^{-\nu-\gamma} \, \|Q(a)\|_{\mathcal{C}_{3+\nu}^\gamma(M)} \\ 
&\leq C \, \varepsilon^{-\nu-\gamma} \, \|D_A^* F_A\|_{\mathcal{C}_{3+\nu}^\gamma(M)} + C \, 
\varepsilon^{-2\nu-\gamma} \, \|a\|_{\mathcal{C}_{1+\nu}^{2,\gamma}(M)}^2 + C \, \varepsilon
^{-3\nu-\gamma} \, \|a\|_{\mathcal{C}_{1+\nu}^{2,\gamma}(M)}^3 \\ 
&\leq C \, \varepsilon^{2-\nu-\gamma} 
\end{align*} 
for all $a \in \mathcal{C}_{1+\nu}^{2,\gamma}(M)$ satisfying 
\[\|a\|_{\mathcal{C}_{1+\nu}^{2,\gamma}(M)} \leq \varepsilon^{\frac{3}{2}}.\] 
Moreover, we have 
\begin{align*} 
\|\Phi(a) - \Phi(a')\|_{\mathcal{C}_{1+\nu}^{2,\gamma}(M)} 
&\leq C \, \varepsilon^{-\nu-\gamma} \, \|Q(a) - Q(a')\|_{\mathcal{C}_{3+\nu}^\gamma(M)} \\ 
&\leq C \, \varepsilon^{\frac{3}{2}-2\nu-\gamma} \, \|a - a'\|_{\mathcal{C}_{1+\nu}^{2,\gamma}
(M)} 
\end{align*} 
for all $a,a' \in \mathcal{C}_{1+\nu}^{2,\gamma}(M)$ satisfying 
\[\|a\|_{\mathcal{C}_{1+\nu}^{2,\gamma}(M)}, \, 
\|a'\|_{\mathcal{C}_{1+\nu}^{2,\gamma}(M)} \leq \varepsilon^{\frac{3}{2}}.\] 
Hence, it follows from the contraction mapping principle that there exists a $1$-form $a \in 
\mathcal{C}_{1+\nu}^{2,\gamma}(M)$ such that 
\[\|a\|_{\mathcal{C}_{1+\nu}^{2,\gamma}(M)} \leq C \, \varepsilon^{2-\nu-\gamma}\] 
and 
\[\Phi(a) = a.\] 
From this it follows that 
\[\mathbb{G} \, (\Id - \mathbb{P}) \, (D_A^* F_A) + a + \mathbb{G} \, (\Id - \mathbb{P}) \, Q(a) 
= 0,\] 
hence 
\[(\Id - \mathbb{P}) \, (D_A^* F_A) + (\Id - \mathbb{P}) \, \mathbb{L}_A a + (\Id - \mathbb{P}) 
\, Q(a) = 0.\] 
Thus, we conclude that 
\[(\Id - \mathbb{P}) \, (D_{\tilde{A}}^* F_{\tilde{A}} + D_{\tilde{A}} D_{\tilde{A}}^* a) = 0.\] 
This proves the assertion. \\

\begin{corollary}
If $\tilde{A}$ satisfies 
\[\mathbb{P} \, (D_{\tilde{A}}^* F_{\tilde{A}} + D_{\tilde{A}} D_{\tilde{A}}^* a) = 0,\] 
then $\tilde{A}$ is a solution of the Yang-Mills equations, i.e. 
\[D_{\tilde{A}}^* F_{\tilde{A}} = 0.\]
\end{corollary}

\textbf{Proof.} 
By definition of $\tilde{A}$, we have 
\[(\Id - \mathbb{P}) \, (D_{\tilde{A}}^* F_{\tilde{A}} + D_{\tilde{A}} D_{\tilde{A}}^* a) = 0.\] 
Hence, if $\tilde{A}$ satisfies 
\[\mathbb{P} \, (D_{\tilde{A}}^* F_{\tilde{A}} + D_{\tilde{A}} D_{\tilde{A}}^* a) = 0,\] 
then we obtain 
\[D_{\tilde{A}}^* F_{\tilde{A}} + D_{\tilde{A}} D_{\tilde{A}}^* a = 0.\] 
Using the Bianchi identity, we obtain 
\[D_{\tilde{A}}^* D_{\tilde{A}} D_{\tilde{A}}^* a = D_{\tilde{A}}^* D_{\tilde{A}}^* F_{\tilde
{A}} + D_{\tilde{A}}^* D_{\tilde{A}} D_{\tilde{A}}^* a = 0.\] 
Integrating over $M$, we conclude that 
\[D_{\tilde{A}} D_{\tilde{A}}^* a = 0,\] 
hence 
\[D_{\tilde{A}}^* F_{\tilde{A}} = 0.\] 
This proves the assertion. \\

\section{The balancing condition}

\begin{proposition}
Let $g_0$ be the product metric on the normal bundle $NS$ (cf. Section 4). Then the fibrewise 
projection $\Pi(D_A^{*_{g_0}} F_A)$ is given by 
\[\Pi(D_A^{*_{g_0}} F_A) = \Big ( \Delta v,\frac{1}{\lambda} \, \Delta \lambda - \frac{1}{4} \, 
|\theta|^2, \frac{1}{\lambda^2} \, \sum_{i=1}^{n-4} \nabla_i(\lambda^2 \, \theta_i) \Big ).\] 
\end{proposition}

\textbf{Proof.} 
Using the results from Section 4, we obtain 
\[\sum_{\beta=1}^4 D_{A,e_\beta^\perp} F_A(e_i,e_\beta^\perp) = 0.\] 
Moreover, the Bianchi identity implies that 
\[D_{A,e_\beta^\perp} F_A(e_i,e_\alpha^\perp) - D_{A,e_\alpha^\perp} F_A(e_i,e_\beta^\perp) + 
D_{A,e_i} F_A(e_\alpha^\perp,e_\beta^\perp) = 0.\] 
From this it follows that 
\begin{align*} 
&\sum_{i=1}^{n-4} \sum_{\alpha,\beta=1}^4 \big ( \nabla_{e_\beta^\perp} \langle F_A(e_i,e_\beta
^\perp),F_A(e_i,e_\alpha^\perp) \rangle - \frac{1}{2} \, \nabla_{e_\alpha^\perp} \langle F_A(e
_i,e_\beta^\perp),F_A(e_i,e_\beta^\perp) \rangle \\ 
&+ \nabla_{e_i} \langle F_A(e_i,e_\beta^\perp),F_A(e_\alpha^\perp,e_\beta^\perp) \rangle \big ) 
\, X^\alpha \\ 
&= \sum_{i=1}^{n-4} \sum_{\alpha,\beta=1}^4 \langle D_{A,e_i} F_A(e_i,e_\beta^\perp),F_A(e
_\alpha^\perp,e_\beta^\perp) \rangle \, X^\alpha \\ 
&= \langle D_A^{*_{g_0}} F_A,F_A(X,\cdot) \rangle. 
\end{align*} 
If $X$ is a vector field of the form 
\[X = \varepsilon \, w_\rho \, e_\rho^\perp + \mu \, (y - \varepsilon v)_\rho \, e_\rho^\perp 
+ r_{\rho\sigma} \, (y - \varepsilon v)_\sigma \, e_\rho^\perp,\] 
then we have 
\begin{align*} 
&\int_{NS_x} \sum_{i=1}^{n-4} \sum_{\alpha,\beta=1}^4 \big ( \nabla_{e_\beta^\perp} \langle F_A
(e_i,e_\beta^\perp),F_A(e_i,e_\alpha^\perp) \rangle - \frac{1}{2} \, \nabla_{e_\alpha^\perp} 
\langle F_A(e_i,e_\beta^\perp),F_A(e_i,e_\beta^\perp) \rangle \big ) \, X^\alpha \\ 
&= -\int_{NS_x} \sum_{i=1}^{n-4} \sum_{\alpha,\beta=1}^4 \big ( \langle F_A(e_i,e_\beta^\perp),
F_A(e_i,e_\alpha^\perp) \rangle \, \nabla_\beta X^\alpha - \frac{1}{2} \, \langle F_A(e_i,e
_\beta^\perp),F_A(e_i,e_\beta^\perp) \rangle \, \nabla_\alpha X^\alpha \big ) \\ 
&= \varepsilon^2 \, (4\pi^2 \, |\nabla v|^2 + 8\pi^2 \, |\nabla \lambda|^2 + 2\pi^2 \, |\theta|
^2) \, \mu. 
\end{align*} 
Similarly, we obtain 
\begin{align*} 
&\int_{NS_x} \sum_{\alpha,\beta=1}^4 \langle F_A(e_i,e_\beta^\perp),F_A(e_\alpha^\perp,e_\beta
^\perp) \rangle \, X^\alpha \\ 
&= -\varepsilon^2 \, (4\pi^2 \, \langle \nabla_i v,w \rangle + 8\pi^2 \, \lambda \, \nabla_i 
\lambda \, \mu + 2\pi^2 \, \lambda^2 \, \langle \theta_i,r \rangle). 
\end{align*} 
Differentiating this identity, we obtain 
\begin{align*} 
&\int_{NS_x} \sum_{i=1}^{n-4} \sum_{\alpha,\beta=1}^4 \nabla_{e_i} \langle F_A(e_i,e_\beta
^\perp),F_A(e_\alpha^\perp,e_\beta^\perp) \rangle \big ) \, X^\alpha \\ 
&= \sum_{i=1}^{n-4} \nabla_{e_i} \int_{NS_x} \sum_{\alpha,\beta=1}^4 \langle F_A(e_i,e_\beta
^\perp),F_A(e_\alpha^\perp,e_\beta^\perp) \rangle \, X^\alpha \\ 
&- \int_{NS_x} \sum_{i=1}^{n-4} \sum_{\alpha,\beta=1}^4 \langle F_A(e_i,e_\beta^\perp),F_A(e
_\alpha^\perp,e_\beta^\perp) \rangle \, \nabla_i X^\alpha \\ 
&= -\varepsilon^2 \, \bigg ( 4\pi^2 \, \langle \Delta v,w \rangle + 8\pi^2 \, \lambda \, \Delta 
\lambda \, \mu + 8\pi^2 \, |\nabla \lambda|^2 \, \mu + 4\pi^2 \, |\nabla v|^2 \, \mu + 2\pi^2 \, 
\Big \langle \sum_{i=1}^{n-4} \nabla_i(\lambda^2 \, \theta_i),r \Big \rangle \bigg ). 
\end{align*} 
Thus, we conclude that 
\begin{align*} 
&\int_{NS_x} \langle D_A^{*_{g_0}} F_A,F_A(X,\cdot) \rangle \\ 
&= -\varepsilon^2 \, \bigg ( 4\pi^2 \, \langle \Delta v,w \rangle + 8\pi^2 \, \lambda \, \langle 
\Delta \lambda,\mu \rangle - 2\pi^2 \, \lambda^2 \, |\theta|^2 \, \mu + 2\pi^2 \, \Big \langle 
\sum_{i=1}^{n-4} \nabla_i(\lambda^2 \, \theta_i),r \Big \rangle \bigg ). 
\end{align*} 
From this the assertion follows. \\

\begin{proposition}
The fibrewise projection $\Pi(D_A^* F_A)$ satisfies the estimate 
\begin{align*} 
\bigg \| &\Pi(D_A^* F_A) \\ 
&- \bigg ( \Delta v_\rho + \sum_{i,j=1}^{n-4} \sum_{\rho,\sigma=1}^4 h_{ij,\rho} \, h_{ij,
\sigma} \, v_\sigma + \sum_{i=1}^{n-4} \sum_{\rho,\sigma=1}^4 R_{i\rho\sigma i} \, v_\sigma, \\ 
&\hspace{7.5mm} \frac{1}{\lambda} \, \Delta \lambda + \frac{1}{4} \sum_{i,j=1}^{n-4} \sum_{\rho
=1}^4 h_{ij,\rho} \, h_{ij,\rho} + \frac{1}{4} \sum_{i=1}^{n-4} \sum_{\rho=1}^4 R_{i\rho\rho i} 
- \frac{1}{4} \, |\theta|^2, \\ 
&\hspace{7.5mm} \frac{1}{\lambda^2} \, \sum_{i=1}^{n-4} \nabla_i (\lambda^2 \, \theta_{i,\rho
\sigma}) \bigg ) \bigg \|_{\mathcal{C}^\gamma(S)} \leq C \, \varepsilon. 
\end{align*} 
\end{proposition}

\textbf{Proof.}
The Riemannian metric satisfies the asymptotic expansion of the form 
\begin{align*} 
g(e_i,e_j) &= \delta_{ij} + 2 \sum_{\rho=1}^4 h_{ij,\rho} \, y_\rho \\ 
&+ \sum_{k=1}^{n-4} \sum_{\rho,\sigma=1}^4 h_{ik,\rho} \, h_{jk,\sigma} \, y_\rho \, y_\sigma \\ 
&- \sum_{\rho,\sigma=1}^4 R_{i\rho\sigma j} \, y_\rho \, y_\sigma + O(|y|^3) \\ 
g(e_i,e_\alpha^\perp) &= O(|y|^2) \\ 
g(e_\alpha^\perp,e_\beta^\perp) &= \delta_{\alpha\beta} - \frac{1}{3} \sum_{\rho,\sigma=1}^4 
R_{\alpha\rho\sigma\beta} \, y_\rho \, y_\sigma + O(|y|^3). 
\end{align*} 
Using this asymptotic expansion, Proposition 7.2 can be deduced from Proposition 7.1. The 
details are left to the reader. \\

\begin{proposition}
The fibrewise projection of $D_{\tilde{A}}^* F_{\tilde{A}} + D_{\tilde{A}} D_{\tilde{A}}^* a$ 
satisfies 
\begin{align*} 
\bigg \| &\Pi(D_{\tilde{A}}^* F_{\tilde{A}} + D_{\tilde{A}} D_{\tilde{A}}^* a) \\ 
&- \bigg ( \Delta v_\rho + \sum_{i,j=1}^{n-4} \sum_{\rho,\sigma=1}^4 h_{ij,\rho} \, h_{ij,
\sigma} \, v_\sigma + \sum_{i=1}^{n-4} \sum_{\rho,\sigma=1}^4 R_{i\rho\sigma i} \, v_\sigma, \\ 
&\hspace{7.5mm} \frac{1}{\lambda} \, \Delta \lambda + \frac{1}{4} \sum_{i,j=1}^{n-4} \sum_{\rho
=1}^4 h_{ij,\rho} \, h_{ij,\rho} + \frac{1}{4} \sum_{i=1}^{n-4} \sum_{\rho=1}^4 R_{i\rho\rho i} 
- \frac{1}{4} \, |\theta|^2, \\ 
&\hspace{7.5mm} \frac{1}{\lambda^2} \, \sum_{i=1}^{n-4} \nabla_i(\lambda^2 \, \theta_{i,\rho
\sigma}) \bigg ) \bigg \|_{\mathcal{C}^\gamma(S)} \leq C \, \varepsilon^{\frac{1}{32}}. 
\end{align*} 
\end{proposition}

\textbf{Proof.} 
Using the estimate 
\[\|a\|_{C_{1+\nu}^{2,\gamma}(M)} \leq C \, \varepsilon^{2-\nu-\gamma},\] 
we obtain 
\begin{align*} 
\|\Pi \, Q(a)\|_{\mathcal{C}^\gamma(S)} 
&\leq C \, \varepsilon^{-2-\nu-\gamma} \, \|Q(a)\|_{\mathcal{C}_{3+\nu}^\gamma(M)} \\ 
&\leq C \, \varepsilon^{-2-2\nu-\gamma} \, \|a\|_{\mathcal{C}_{1+\nu}^{2,\gamma}(M)}^2 + C \, 
\varepsilon^{-2-3\nu-\gamma} \, \|a\|_{\mathcal{C}_{1+\nu}^{2,\gamma}(M)}^3 \\ 
&\leq C \, \varepsilon^{2-4\nu-3\gamma}. 
\end{align*} 
Moreover, we have 
\[\|\Pi \mathbb{L}_A a\|_{\mathcal{C}^\gamma(S)} \leq C \, \varepsilon^{\frac{1}{32}}.\] 
Hence, the assertion follows from Proposition 7.2. \\

\textbf{Proof of Theorem 1.1.} 
Let 
\[\Xi_\varepsilon(v,\lambda,J,\omega) = \Pi(D_{\tilde{A}}^* F_{\tilde{A}} + D_{\tilde{A}} 
D_{\tilde{A}}^* a).\] 
The first part of Theorem 1.1 follows from Corollary 6.2, the second part from Proposition 7.3. 
\\

\section{An example}

Suppose that the normal bundle $NS$ can be endowed with a $SU(2)$-structure $(J,\omega)$ which 
is parallel with respect to the Levi-Civita connection $\nabla$. This implies that $\theta = 0$. 
Moreover, suppose that $\lambda$ is a postive function on $S$ which satisfies the linear PDE 
\[\Delta \lambda + \frac{1}{4} \sum_{i,j=1}^{n-4} \sum_{\rho=1}^4 h_{ij,\rho} \, h_{ij,\rho} \, 
\lambda + \frac{1}{4} \sum_{i=1}^{n-4} \sum_{\rho=1}^4 R_{i\rho\rho i} \, \lambda = 0.\] 
In addition, we assume that the following non-degeneracy conditions hold: \\

(i) The Jacobi operator of $S$ is invertible. \\

(ii) The kernel of the operator 
\[\Delta + \frac{1}{4} \sum_{i,j=1}^{n-4} \sum_{\rho=1}^4 h_{ij,\rho} \, h_{ij,\rho} + \frac
{1}{4} \sum_{i=1}^{n-4} \sum_{\rho=1}^4 R_{i\rho\rho i}\] 
is spanned by the function $\lambda$. \\

\begin{proposition}
Let $(w,\mu,r)$ be a section of the vector bundle $NS \oplus \mathbb{R} \oplus \Lambda_+^2 NS$ 
such that 
\[\|w\|_{\mathcal{C}^{2,\gamma}(S)} \leq \varepsilon^{\frac{1}{64}},\] 
\[\|\mu\|_{\mathcal{C}^{2,\gamma}(S)} \leq \varepsilon^{\frac{1}{64}},\] 
\[\|s\|_{\mathcal{C}^{2,\gamma}(S)} \leq \varepsilon^{\frac{1}{64}}.\] 
Then the connection $\tilde{A} = A + a$ corresponding to $(w,\mu,r)$ satisfies the estimate 
\begin{align*} 
\bigg \| &\Pi(D_{\tilde{A}}^* F_{\tilde{A}} + D_{\tilde{A}} D_{\tilde{A}}^* a) \\ 
&- \bigg ( \Delta w_\rho + \sum_{i,j=1}^{n-4} \sum_{\rho,\sigma=1}^4 h_{ij,\rho} \, h_{ij,
\sigma} \, w_\sigma + \sum_{i=1}^{n-4} \sum_{\rho,\sigma=1}^4 R_{i\rho\sigma i} \, w_\sigma, \\ 
&\hspace{7.5mm} \frac{1}{\lambda} \, \Delta (\lambda \, \mu) + \frac{1}{4} \sum_{i,j=1}^{n-4} 
\sum_{\rho=1}^4 h_{ij,\rho} \, h_{ij,\rho} \, \mu + \frac{1}{4} \sum_{i=1}^{n-4} \sum_{\rho=1}^4 
R_{i\rho\rho i} \, \mu, \\ 
&\hspace{7.5mm} \frac{1}{\lambda^2} \, \sum_{i=1}^{n-4} \nabla_i (\lambda^2 \, \nabla_i r) \bigg 
) \bigg \|_{\mathcal{C}^\gamma(S)} \leq C \, \varepsilon^{\frac{1}{32}}. 
\end{align*} 
\end{proposition}

\textbf{Proof.} 
This follows immediately from Proposition 7.3. \\

For abbreviation, let 
\begin{align*} 
J(w,\mu,r) = \bigg ( 
&\Delta w_\rho + \sum_{i,j=1}^{n-4} \sum_{\rho,\sigma=1}^4 h_{ij,\rho} \, h_{ij,\sigma} \, 
w_\sigma + \sum_{i=1}^{n-4} \sum_{\rho,\sigma=1}^4 R_{i\rho\sigma i} \, w_\sigma, \\ 
&\frac{1}{\lambda} \, \Delta (\lambda \, \mu) + \frac{1}{4} \sum_{i,j=1}^{n-4} \sum_{\rho=1}^4 
h_{ij,\rho} \, h_{ij,\rho} \, \mu + \frac{1}{4} \sum_{i=1}^{n-4} \sum_{\rho=1}^4 R_{i\rho\rho i} 
\, \mu, \\ 
&\frac{1}{\lambda^2} \, \sum_{i=1}^{n-4} \nabla_i (\lambda^2 \, \nabla_i r) \bigg ). 
\end{align*} 
If $(w,\mu,r)$ is a section of the vector bundle $NS \oplus \mathbb{R} \oplus \Lambda_+^2 NS$ of 
class $\mathcal{C}^{2,\gamma}$, then $J(w,\mu,r)$ is a section of the vector bundle $NS \oplus 
\mathbb{R} \oplus \Lambda_+^2 NS$ of class $\mathcal{C}^\gamma$. \\

The kernel of $J$ is a vector space of dimension $4$. It consists of all triplets $(0,\mu,r)$, 
where $\mu$ is constant on $S$ and $r$ is a parallel section of the vector bundle $\Lambda_+^2 
NS$. \\

\begin{proposition}
There exists a section $(w,\mu,r)$ of the vector bundle $NS \oplus \mathbb{R} \oplus \Lambda_+^2 
NS$ such that 
\[\Pi(D_{\tilde{A}}^* F_{\tilde{A}} + D_{\tilde{A}} D_{\tilde{A}}^* a) \in V,\] 
where $\tilde{A} = A + a$ is the connection corresponding to $(w,\mu,r)$.
\end{proposition}

\textbf{Proof.} 
By Proposition 8.1, we may write 
\[\Pi(D_{\tilde{A}}^* F_{\tilde{A}} + D_{\tilde{A}} D_{\tilde{A}}^* a) = J(w,\mu,r) + R(w,\mu,
r),\] 
where $\|R(w,\mu,r)\|_{\mathcal{C}^\gamma(S)} \leq C \, \varepsilon^{\frac{1}{32}}$ for $\|(w,
\mu,r)\|_{\mathcal{C}^{2,\gamma}(S)} \leq \varepsilon^{\frac{1}{64}}$. Hence, the operator $-J
^{-1} \, R$ maps a ball of radius $\varepsilon^{\frac{1}{64}}$ in the Banach space $\mathcal{C}
^{2,\gamma}(S)$ into a ball of radius $C \, \varepsilon^{\frac{1}{32}}$ in $\mathcal{C}^{2,
\gamma}(S)$. Using an appropriate sequence of smoothing operators, we may approximate the 
mapping $-J^{-1} \, R$ by a sequence of compact mappings. Each of these mappings has a fixed 
point in $\mathcal{C}^{2,\gamma}(S)$ by Schauder's fixed point theorem. Taking limits, we obtain 
a fixed point of the original mapping $-J^{-1} \, R$ in the Banach space $\mathcal{C}^{2,\frac
{\gamma}{2}}(S)$. With this choice of the glueing data $(w,\mu,r)$, we obtain $(w,\mu,r) + 
J^{-1} \, R(w,\mu,r) = 0$, hence $J(w,\mu,r) + R(w,\mu,r) \in V$. \\


\begin{thebibliography}{99}
\bibitem{Ba}
A. Bahri, \textit{An invariant for Yamabe-type flows with applications to scalar curvature 
problems in high dimension,} Duke Math. J. 81, 323-466 (1996)

\bibitem{BC}
A. Bahri and J.M. Coron, \textit{On a nonlinear elliptic equation involving the critical Sobolev 
exponent: The effect of the topology of the domain,} Comm. Pure Appl. Math. 41, 253-290 (1988)

\bibitem{DK}
S. K. Donaldson and P. B. Kronheimer, \textit{The Geometry of Four-Manifolds,}
Oxford University Press (1990)

\bibitem{Ka1}
N. Kapouleas, \textit{Complete constant mean curvature surfaces in Euclidean three-space,} Ann. 
of Math. 131, 239-330 (1990)

\bibitem{Ka2}
N. Kapouleas, \textit{Compact constant mean curvature surfaces in Euclidean three-space,} J. 
Diff. Geom. 33, 683-715 (1991)

\bibitem{KMP}
R. Kusner, R. Mazzeo, and D. Pollack, \textit{The moduli space of of complete embedded constant 
mean curvature surfaces,} Geom. Funct. Anal. 6, 120-137 (1996)

\bibitem{Li}
F.-H. Lin, \textit{Complex Ginzburg-Landau equations and dynamics of vortices, filaments, and 
codimension-$2$ submanifolds,} Comm. Pure Appl. Math. 51, 385-441 (1998)

\bibitem{LR}
F.-H. Lin and T. Rivi\`ere, \textit{Complex Ginzburg-Landau equations in high dimensions and 
codimension two area minimizing currents,} J. Eur. Math. Soc. 1, 237-311 (1999)

\bibitem{MP1}
R. Mazzeo and F. Pacard, \textit{A construction of singular solutions for a semilinear elliptic 
equation using asymptotic analysis,} J. Diff. Geom. 44, 331-370 (1996)

\bibitem{MP2}
R. Mazzeo and F. Pacard, \textit{Constant mean curvature surfaces with Delaunay ends,} Comm. 
Anal. Geom. 9, 169-237 (2001)

\bibitem{MPP}
R. Mazzeo, F. Pacard and D. Pollack, \textit{Connected sums of constant mean curvature surfaces 
in Euclidean $3$-space,} J. Reine Angew. Math. 536, 115-165 (2001)

\bibitem{MS}
R. Mazzeo and N. Smale, \textit{Conformally flat metrics of constant positive scalar curvature 
on subdomains of the sphere,} J. Diff. Geom. 34, 581-621 (1991)

\bibitem{PR1}
F. Pacard and M. Ritor\'e, \textit{From constant mean curvature hypersurfaces to the gradient 
theory of phase transitions,} preprint (2003)

\bibitem{PR2}
F. Pacard and T. Rivi\`ere, \textit{Linear and Nonlinear Aspects of Vortices. The 
Ginzburg-Landau Model,} Progress in Nonlinear Differential Equations and their Applications, 
vol. 39, Birkh\"auser, Boston (2000)

\bibitem{Po}
D. Pollack, \textit{Nonuniqueness and high energy solutions for a conformally
invariant scalar equation,} Comm. Anal. Geom. 3, 347-414 (1993)

\bibitem{TT}
T. Tao and G. Tian, \textit{A singularity removal theorem for Yang-Mills fields in higher 
dimensions,} preprint (2002)

\bibitem{Ta1}
C. H. Taubes, \textit{Self-dual Yang-Mills connections on non-self-dual $4$-manifolds,} J. Diff. 
Geom. 17, 139-170 (1982)

\bibitem{Ta2}
C. H. Taubes, \textit{$\text{SW} \Longrightarrow \text{Gr}$: from the Seiberg-Witten equations 
to pseudo-holomorphic curves,} J. Amer. Math. Soc. 9, 845-918 (1996)

\bibitem{Ta3}
C. H. Taubes, \textit{$\text{Gr} \Longrightarrow \text{SW}$: from pseudo-holomorphic curves to 
Seiberg-Witten solutions,} J. Diff. Geom. 51, 203-334 (1999)

\bibitem{Ti}
G. Tian, \textit{Gauge theory and calibrated geometry,} Ann. of Math. 151, 193-268 (2000)
\end{thebibliography}
\end{document}